\documentclass[12pt]{article}
\usepackage{amssymb}
\usepackage{latexsym}

\newtheorem{claim}{Claim}[section]
\newtheorem{conj}{Conjecture}[section]
\newtheorem{cor}{Corollary}[section]
\newtheorem{defn}{Definition}[section]
\newtheorem{hyp}{Hypothesis}[section]
\newtheorem{lemma}{Lemma}[section]
\newtheorem{prop}{Proposition}[section]

\newtheorem{thm}{Theorem}[section]

\newcounter{example}[section]
\newcounter{rem}[section]

\newcommand{\barbbQ}{\ensuremath{\bar{\bbQ}}}
\newcommand{\barbbQp}{\ensuremath{\bar{\bbQ}_p}}

\newcommand{\bbC}{\ensuremath{\mathbb{C}}}

\newcommand{\bbN}{\ensuremath{\mathbb{N}}}

\newcommand{\bbQ}{\ensuremath{\mathbb{Q}}}

\newcommand{\bbR}{\ensuremath{\mathbb{R}}}
\newcommand{\bbZ}{\ensuremath{\mathbb{Z}}}
\newcommand{\bbZp}{\ensuremath{\mathbb{Z}_p}}
\newcommand{\be}{\ensuremath{{\bf e}}}
\newcommand{\beql}[1]{\begin{equation}\label{#1}}
\newcommand{\bPf}{\noindent \textsc{Proof\ }}

\newcommand{\bv}{\ensuremath{{\bf v}}}
\newcommand{\bw}{\ensuremath{{\bf w}}}

\newcommand{\cB}{\ensuremath{\mathcal{B}}}
\newcommand{\cC}{\ensuremath{\mathcal{C}}}
\newcommand{\cD}{\ensuremath{\mathcal{D}}}

\newcommand{\cf}{\emph{cf.}}

\newcommand{\cH}{\ensuremath{\mathcal{H}}}
\newcommand{\cI}{\ensuremath{\mathcal{I}}}
\newcommand{\cJ}{\ensuremath{\mathcal{J}}}

\newcommand{\Cl}{{\rm Cl}}

\newcommand{\cO}{\ensuremath{\mathcal{O}}}
\newcommand{\cP}{\ensuremath{\mathcal{P}}}
\newcommand{\cQ}{\ensuremath{\mathcal{Q}}}
\newcommand{\cR}{\ensuremath{\mathcal{R}}}
\newcommand{\cS}{\ensuremath{\mathcal{S}}}
\newcommand{\cT}{\ensuremath{\mathcal{T}}}

\newcommand{\cV}{\ensuremath{\mathcal{V}}}

\newcommand{\displaymapdef}[5]
{\[
\begin{array}{rcrcl}
 #1 &:& #2 &\longrightarrow& #3 \\
    & &    &                    \\
    & & #4 &\longmapsto    & #5
\end{array}
\]}

\newcommand{\eeq}{\end{equation}}
\newcommand{\eg}{\emph{e.g.}}
\newcommand{\ePf}{\hspace*{\fill}~$\Box$\vertsp\par}

\newcommand{\etc}{\emph{etc.}}
\newcommand{\example}{\refstepcounter{example}
\noindent{\sc Example \theexample}}
\newcommand{\fa}{\ensuremath{\mathfrak{a}}}

\newcommand{\fc}{\ensuremath{\mathfrak{c}}}
\newcommand{\ff}{\ensuremath{\mathfrak{f}}}
\newcommand{\fg}{\ensuremath{\mathfrak{g}}}

\newcommand{\fm}{\ensuremath{\mathfrak{m}}}
\newcommand{\fn}{\ensuremath{\mathfrak{n}}}
\newcommand{\fp}{\ensuremath{\mathfrak{p}}}
\newcommand{\fP}{\ensuremath{\mathfrak{P}}}
\newcommand{\fq}{\ensuremath{\mathfrak{q}}}
\newcommand{\fQ}{\ensuremath{\mathfrak{Q}}}

\newcommand{\fs}{\ensuremath{\mathfrak{s}}}
\newcommand{\fS}{\ensuremath{\mathfrak{S}}}

\newcommand{\fw}{\ensuremath{\mathfrak{w}}}
\newcommand{\fW}{\ensuremath{\mathfrak{W}}}
\newcommand{\fz}{\ensuremath{\mathfrak{z}}}
\newcommand{\Gal}{{\rm Gal}}

\newcommand{\half}{\frac{1}{2}}

\newcommand{\ie}{\emph{i.e.}}
\newcommand{\inv}{^{-1}}
\newcommand{\ndiv}{\nmid}
\newcommand{\nin}{\not\in}
\newcommand{\ord}{{\rm ord}}

\newcommand{\refeq}[1]{~(\ref{eq:#1})}
\newcommand{\rem}{\refstepcounter{rem}\noindent{\sc Remark \therem}}

\renewcommand{\theexample}{\thesection.\arabic{example}}
\renewcommand{\therem}{\thesection.\arabic{rem}}
\newcommand{\Tr}{{\rm Tr}}
\newcommand{\vertsp}{\vspace{1ex}}

\newcommand{\Clmk}{\ensuremath{\Cl_\fm(k)}}
\newcommand{\Hrhov}{H($\rho$,$\bv$)}

\newcommand{\ktuinf}{\ensuremath{k^\times_{\uinf}}}
\newcommand{\Mmi}{\ensuremath{M_{\fm,i}}}

\newcommand{\PhiKk}{\ensuremath{\Phi_{K/k}}}
\newcommand{\PhiKko}{\ensuremath{\Phi_{K/k}(0)}}
\newcommand{\PhiKkoa}{\ensuremath{\Phi_{K/k}(0)^\ast}}
\newcommand{\PhiKks}{\ensuremath{\Phi_{K/k}(s)}}
\newcommand{\Phim}{\ensuremath{\Phi_{\fm}}}
\newcommand{\pnpo}{\ensuremath{{p^{n+1}}}}
\newcommand{\sKk}{\ensuremath{\fs_{K/k}}}
\newcommand{\SKk}{\ensuremath{\fS_{K/k}}}

\newcommand{\tff}{\ensuremath{\tilde{\ff}}}

\newcommand{\tfm}{\ensuremath{{\tilde{\fm}}}}
\newcommand{\tfmi}{\ensuremath{{\tfm^{(i)}}}}
\newcommand{\Tmi}{\ensuremath{T_{\fm,i}}}
\newcommand{\Ttmi}{\ensuremath{T_{\tfm,i}}}

\newcommand{\tK}{\ensuremath{\tilde{K}}}
\newcommand{\tphi}{\ensuremath{\tilde{\phi}}}
\newcommand{\ThetaKk}{\ensuremath{\Theta_{K/k}}}
\newcommand{\tsigma}{\ensuremath{\tilde{\sigma}}}
\newcommand{\uinf}{\ensuremath{{\underline{\infty}}}}
\newcommand{\ujct}{\ensuremath{\underline{j\tau}}}
\newcommand{\ul}{\ensuremath{{\underline{l}}}}
\newcommand{\um}{\ensuremath{{\underline{m}}}}
\newcommand{\un}{\ensuremath{{\underline{n}}}}
\newcommand{\uo}{\ensuremath{{\underline{1}}}}
\newcommand{\ur}{\ensuremath{{\underline{r}}}}
\newcommand{\us}{\ensuremath{{\underline{s}}}}
\newcommand{\ut}{\ensuremath{{\underline{t}}}}
\newcommand{\utau}{\ensuremath{{\underline{\tau}}}}
\newcommand{\uX}{\ensuremath{{\underline{X}}}}
\newcommand{\uz}{\ensuremath{{\underline{z}}}}
\newcommand{\uzero}{\ensuremath{{\underline{0}}}}
\newcommand{\uzeta}{\ensuremath{{\underline{\zeta}}}}
\begin{document}
\title{On Twisted Zeta-Functions at $s=0$}
\author{D. Solomon
\\King's College London\\
{\tt solomon@mth.kcl.ac.uk}}
\maketitle
\section{Introduction}
Let $K/k$ be an abelian extension of number fields with Galois group $G$ and
for each $\sigma\in G$ let $\zeta(s;\sigma)$ denote the associated partial zeta-function.
Letting $\sigma$ vary, we obtain a meromorphic, group-ring-valued function
$\Theta_{K/k}:\bbC\rightarrow \bbC G$ by setting
$\Theta_{K/k}(s):=\sum_{\sigma\in G}\zeta(s;\sigma)\sigma\inv$
(for more details of this construction, see Section~\ref{sec:thetaandphi}).
We first note the \emph{rationality property} satisfied by the value at $s=0$, namely
$\Theta_{K/k}(0)\in\bbQ G$. This
follows from work of Siegel~\cite{Si} and Klingen (or of
Shintani~\cite[Cor. to Thm. 1]{Shi}).
Next, let $\mu(K)$ denote the group of roots of unity in $K$
and let ${\rm ann}_{\bbZ G}(\mu(K))$ be its annihilator ideal as a module for the
group-ring $\bbZ G$.
Deligne-Ribet and Pi.~Cassou-Nogu\`es proved the following
\emph{integrality property} concerning
$\Theta_{K/k}(0)$ (see~\cite{D-R} and~\cite{CN}):
\beql{eq:introA}
{\rm ann}_{\bbZ G}(\mu(K))\ThetaKk(0)\subset\bbZ G
\eeq
In particular, if $w_K=|\mu(K)|$ then  $w_K\ThetaKk(0)$ lies in $\bbZ G$.
One can show that $\Theta_{K/k}(0)$ is essentially always zero unless
$k$ is totally real and $K$ is totally complex and that in this case  $(1+c)\Theta_{K/k}(0)$
vanishes for every complex conjugation $c\in G$.
(We say that $\Theta_{K/k}(0)$ lies in the `minus part' of $\bbQ G$ w.r.t.\ such $c$).
In fact, one loses little generality in assuming that $K$ is a CM field with
whose unique complex conjugation $c$ lies in $G$ (see Remark~\ref{rem:2A.5}(ii)).
In this context,
the \emph{Brumer-Stark Conjecture} takes the form of
a conjectural generalisation of Stickelberger's Theorem
and may be stated as follows. If $\fa$ is any fractional ideal of $K$, then
$\fa^{w_K\ThetaKk(0)}$ is a principal ideal generated by some element $a$ in the minus part
of $K^\times$ w.r.t.\ $c$.
(One also imposes other conditions on $a$. See~\cite{Gr} for more details
and an account of recent work demonstrating this conjecture in many cases.)

This paper initiates the study at $s=0$ of a function
$\Phi_{K/k}(s)$ related to $\Theta_{K/k}(s)$ and in which, roughly speaking,
the partial zeta-functions are replaced by the \emph{twisted
zeta-functions} introduced and studied in~\cite{twizas} --~\cite{sconzp} and~\cite{zetap2}
(see \emph{ibid.} and Section~\ref{sec:thetaandphi}).
We shall establish what, in the broadest terms, may be described as
a \emph{`rationality property'} and an \emph{`integrality property'}
for $\PhiKko$ and also formulate a new \emph{`Stark-type conjecture'} related to it.
The statements of these properties and our conjecture
are, however, far from being obvious variants of the corresponding statements
for $\Theta_{K/k}(0)$ mentioned above. In particular
they are all fundamentally $p$-adic in nature.

Here is a slightly more detailed summary of this paper.
Section~\ref{sec:thetaandphi} contains more on the definitions
of the functions $\ThetaKk(s)$ and $\PhiKk(s)$ and gives
two relations between them. The first is
a `functional equation' relating $\ThetaKk(s)$ to $\PhiKk(1-s)$ (Theorem~\ref{thm:2A}).
It means that one could equally well regard this paper as initiating the study
of the function $\ThetaKk(s)$ at $s=1$ (but see Remark~\ref{rem:2A.5}(iv)).
Our second relation, Theorem~\ref{thm:2B}, is, roughly speaking, a formula
linking $\Phi_{K/k}(s)$ at to the functions $\Theta_{\tilde{K}/k}(s)$
as $\tilde{K}/k$ runs through certain sub-extensions of $K/k$.

As in the case of $\ThetaKk(0)$, our study of $\PhiKko$ loses
very little from the assumption that $k$ is totally real which we shall
therefore make for the rest of this summary.
The `rationality property' mentioned above is Proposition~\ref{prop:2D}.
In order to state it, we first define a $p$-adic
group-ring-valued regulator $R_{K/k,p}$ on the group
$U_p(K)$ of $p$-semilocal units of $K$. Then we show that
the product $\sqrt{d_k}\PhiKk(0)^\ast R_{K/k,p}(\theta)$ lies in $\bbQ_p G$  for any
$\theta\in\bigwedge^d_{\bbZ G}U_p(K)$, where $d_k$ denotes the absolute discriminant of $k$
and other notations will be defined later.
If we denote this product by $\sKk(\theta)$ then we have defined for each prime $p$ a map
$\sKk$ from $\bigwedge^d_{\bbZ G}U_p(K)$ to $\bbQ_p G$.
The image $\SKk$ of $\sKk$ is a $\bbZ_p G$ submodule of $\bbQ_p G$ which is
somewhat analogous to the generalised Stickelberger ideal
${\rm ann}_{\bbZ G}(\mu(K))\ThetaKk(0)$.
Unlike $\Theta_{K/k}(s)$ however, $\PhiKk(s)$ has no `trivial zeroes' at $s=0$. Consequently
$\SKk$ has the pleasing property of spanning the entire minus part of $\bbQ_p G$
(w.r.t\ all complex conjugations.
See Remark~\ref{rem:another}(ii) for more details).

Section~\ref{sec:integrality} contains our main result.
 This is Theorem~\ref{thm:3A} -- the
`integrality property' for $\SKk$ -- which implies in particular that
$\SKk$ is contained in $\bbZ_p G$, provided that certain hypotheses are satisfied
(principally that $p$ is odd and splits completely in $k$). It is not yet clear to
what extent these hypotheses are necessary for the conclusion of Theorem~\ref{thm:3A}. They
do however figure prominently in
its proof which is by far the most substantial in this paper and
draws on two different sources.
On the one hand it borrows ideas from Coleman's method in~\cite{Coleman1} for studying
the the dual of the image of the local logarithm
using one-variable $p$-adic formal power-series.
On the other, it uses Shintani's method (with improvements from Colmez) for generating
twisted zeta-values by means of cone decompositions and
multivariable formal power series. The algebraic properties of the
two sets of power-series marry together very naturally and actually lead to
a neat formula for $\sKk(\theta)$ under our hypotheses.

The final section of this paper contains further discussion of integrality questions
and three conjectures. The last of these, Conjecture~\ref{conj:4D}, is `of Stark-type' and
was motivated by the results of~\cite{sconzp}.
Assuming that $p\neq 2$ splits in $k$, that $K$ contains the $\pnpo$th roots of unity
and that $\SKk\subset\bbZ_p G$, it
proposes certain congruences for $\sKk(\theta)$ modulo $\pnpo$ in terms of Hilbert symbols
and the Rubin-Stark units of $K^+/k$, where $K^+$ denotes the maximal
real subfield of $K$.

In addition to those introduced above, the following
notations and conventions will be used throughout this paper.
All number fields will be considered
as finite extensions of $\bbQ$ within the algebraic closure $\barbbQ$ of $\bbQ$
in $\bbC$. Concerning the `base field' $k$, we shall write
$\cO=\cO_k$ for its ring of integers, $E(k)=\cO^\times$ for its unit group and
$S_\infty=S_\infty(k)$ for the set of its
infinite places, $r_1(k)$ of which we shall initially assume to be real and $r_2(k)$ complex
so that $r_1(k)+2r_2(k)=d:=[k:\bbQ]$. We also fix once and for all elements
$\tau_1,\ldots,\tau_d$ of $\Gal(\barbbQ/\bbQ)$ extending the $d$ distinct embeddings
$k\rightarrow\barbbQ$.

A \emph{cycle} for $k$ will be a formal product $\ff\fz$
where $\ff$ is a non-zero ideal of $\cO$ and $\fz$ is (the formal
product of) a subset of the real places of $k$. The ray-class
group and the ray-class field of $k$ modulo $\fm$ will be denoted
$\Cl_\fm(k)$
and $k(\fm)$ respectively. If $G_\fm$  denotes $\Gal(k(\fm)/k)$ then the Artin isomorphism
$\Cl_\fm(k)\rightarrow G_\fm$ takes a ray-class $\fc$ to $\sigma_{\fc}=\sigma_{\fc,\fm}$.
If $\fa$ is a fractional ideal of $k$ prime to $\ff$ then its class in $\Cl_\fm(k)$
will be denoted $[\fa]_\fm$ and we shall sometimes write $\sigma_{\fa,\fm}$ in place of
$\sigma_{[\fa]_\fm,\fm}$. If $\tfm$ is a cycle dividing $\fm$, then $k(\tfm)$ is contained in
$k(\fm)$ and the restriction homomorphism
$G_\fm\rightarrow G_\tfm$ corresponds by the Artin maps to the homomorphism
$\Cl_{\fm}(k)\rightarrow\Cl_\tfm(k)$ taking $[\fa]_\fm$ to
$[\fa]_\tfm$. We shall write $\pi_{\fm,\tfm}$ for either homomorphism or
indeed for the $R$-linear extension of $\pi_{\fm,\tfm}$ to a homomorphism of group rings
$R G_\fm\rightarrow R G_\tfm$ for any commutative ring $R$.

If $K$ is any abelian extension of $k$ we shall, by a slight abuse of notation,
write $S_{\rm ram}(K/k)$ for the set of finite places of
$k$ which ramify in $K$ together with \emph{all those in $S_\infty$}.
For any place $v$ of $k$, whether finite or infinite,
$D_v(K/k)$ will denote the decomposition subgroup of $G$
associated to some (hence any) place of $K$ above $v$.
If $\tilde{K}/k$ is a subextension of $K/k$ then $\pi_{K,\tilde{K}}$
will denote the restriction $\Gal(K/k)\rightarrow \Gal(\tilde{K}/k)$
(linearly extended to group rings where appropriate).
We shall write $\fm(K)=\ff(K)\fz(K)$ for the
\emph{conductor} of $K$ over $k$, namely the minimal cycle $\fm$ such that $K\subset k(\fm)$.
The support of $\fm(K)$ consists of the places ramified in $K/k$ (finite \emph{or} infinite)
and for any fractional ideal
$\fa$ of $k$ whose support is disjoint from this set, we write $\sigma_{\fa,K}$ for
the corresponding element of $G$ under the Artin map, \ie\
$\sigma_{\fa,K}=\pi_{k(\fm(K)),K}(\sigma_{\fa,\fm(K)})$.

For each prime number $p$ we denote by
$\bbC_p$ the completion of a  fixed algebraic closure $\barbbQ_p$ of $\bbQ_p$
with respect to the $p$-adic metric. We denote by $|\cdot|_p$ the unique absolute value
on $\bbC_p$ normalised such that $|p|_p=p\inv$. Finally, for any positive integer
$f$, we shall write $\mu_f$ for the group of $f$th roots of unity, whether in $\bbC$, or in
$\bbC_p$ for some $p$.

It is my pleasure to thank Martin Taylor for a useful discussion during the
preparation of this paper.

\section{The Functions $\Theta_{K/k}(s)$ and $\Phi_{K/k}(s)$}\label{sec:thetaandphi}
We first record some basic facts about these functions, referring
to~\cite[Ch. IV]{Tate} and to~\cite{twizas},~\cite{zetap1},~\cite{sconzp}
for more details on $\Theta_{K/k}(s)$ and $\Phi_{K/k}(s)$ respectively.
If $\fm=\ff\fz$ is any cycle for $k$ and $\fc$ is any ideal class in $\Clmk$, we
define the corresponding partial zeta-function of a complex variable $s$
by the following Dirichlet series (absolutely convergent for ${\rm Re}(s)>1$)
\[
\zeta(s;\fc)=\zeta_\fm(s;\fc):=\sum_{\fa}N\fa^{-s}
\]
where the sum runs over the set of integral ideals $\fa$
(prime to $\ff$) in $\fc$.
This extends to a meromorphic function on $\bbC$ having only a
simple pole at $s=1$ and we define a $\bbC G_\fm$-valued meromorphic function
$\Theta_\fm$ on $\bbC$ by setting
\beql{eq:2Z}
\Theta_\fm(s):=\sum_{\fc\in\Clmk}\zeta(s;\fc)\sigma_{\fc,\fm}\inv=
\prod_{\fp\ndiv\ff}(1-N\fp^{-s}\sigma_{\fp,\fm}\inv)\inv
\eeq
where the Euler product (over primes ideals of $k$ not dividing $\ff$)
converges for ${\rm Re}(s)>1$.
If $\tilde{\fm}=\tilde{\ff}\tilde{\fz}$ is another cycle with $\tfm|\fm$ then clearly
\beql{eq:2A}
\pi_{\fm,\tfm}(\Theta_{\fm}(s))=
\prod_{\fp|\ff,\ \fp\ndiv\tff}
(1-N\fp^{-s}\sigma_{\fp,\tfm}\inv)\Theta_{\tfm}(s)
\eeq
For any (finite) abelian extension $K$ of $k$ with group $G$ we set
\beql{eq:2A.25}
\Theta_{K/k}(s):=\pi_{k(\fm(K)),K}(\Theta_{\fm(K)}(s))
\eeq
(a meromorphic function of $s\in\bbC$ with values in $\bbC G$). We thus
have the following expressions \emph{intrinsic to K}
\beql{eq:2A.5}
\Theta_{K/k}(s)=
\prod_{\fp\nin S_{\rm ram}(K/k)}(1-N\fp^{-s}\sigma_{\fp,K}\inv)\inv=
\sum_{\sigma\in G}\zeta(s;\sigma)\sigma\inv
\eeq
as in the introduction, where, for ${\rm Re}(s)>1$,
the partial zeta-function $\zeta(s;\sigma)=\zeta_{K/k}(s;\sigma)$
equals the sum $\sum N\fa^{-s}$
as $\fa$ runs through the ideals of $\cO$ with support disjoint
from $S_{\rm ram}(K/k)$ and such that $\sigma_{\fa,K}=\sigma$. If $\tilde{K}/k$ is
a sub-extension of $K/k$ then clearly
\beql{eq:2A.75}
\pi_{K,\tilde{K}}(\Theta_{K/k}(s))=
\prod_{\fp\in S_{\rm ram}(K/k)\atop \fp\nin S_{\rm ram}(\tilde{K}/k)}
(1-N\fp^{-s}\sigma_{\fp,\tilde{K}}\inv)\Theta_{\tilde{K}/k}(s)
\eeq
\rem\ Note that $\Theta_\fm$ and $\Theta_{K/k}$ are the functions
$\Theta_{k(\fm)/k,S_\infty\cup S_\ff}$ and
$\Theta_{K/k, S_{\rm ram}(K/k)}$ of~\cite[Ch. IV]{Tate}, where
$S_\ff:=\{\fp:\fp|\ff\}$.
The containment $S_\ff\supset S_{\rm ram}(k(\fm)/k)$
may be strict for a general cycle $\fm$ (though not if $\fm$ is a conductor)
in which case it follows that
$\Theta_{\fm}(s)$ is not equal to
$\Theta_{k(\fm)/k}(s)$ but is obtained from it by multiplying by
$(1-N\fp^{-s}\sigma_{\fp,k(\fm)}\inv)$ for all
$\fp\in S_\ff\setminus S_{\rm ram}(k(\fm)/k)$.
In this case $\zeta_\fm(s;\fc)$ is not equal to $\zeta_{k(\fm)/k}(s;\sigma_{\fc,\fm})$,
for $\fc\in \Clmk$.\vspace{2ex}\\
For any cycle $\fm=\ff\fz$, we defined in~\cite{twizas} a finite
set denoted $\fW_\fm$ and consisting, in brief, of the
$\fz$-equivalence classes of those additive characters of
fractional ideals of $k$ whose precise $\cO$-annihilator is $\ff$.
We equipped $\fW_\fm$ with a distinguished element, denoted $\fw_\fm^0$,
and a free, transitive action of $\Clmk$. For each
$\fw\in\fW_\fm$ we shall here write simply $Z(s;\fw)$ for the
case $T=\emptyset$ (the empty set) of the twisted zeta-function
$Z_T(s;\fw)$. We refer to~\cite{twizas} for the precise definition of this latter
Dirichlet series (see also the proofs of Theorem~\ref{thm:2B} and
Lemma~\ref{lemma:3K}, and Example~\ref{ex:2A} for a special case). As in~\cite{twizas}, we set
\beql{eq:2A.9}
\Phi_\fm(s):=\sum_{\fc\in\Clmk}Z(s;\fc\cdot\fw_\fm^0)\sigma_\fc\inv
\eeq
to get a meromorphic, $\bbC G_\fm$-valued function $\Phi_{\fm}$
with at most a simple pole at $s=1$. (In fact, $\Phi_{\fm}$ is holomorphic if $\ff\neq \cO$).
For $\tfm|\fm$,
Theorem~3.2 of~\cite{twizas} gives the following analogue of\refeq{2A}
\beql{eq:2B}
\pi_{\fm,\tfm}(\Phi_{\fm}(s))=\left(\frac{N\ff}{N\tff}\right)^{1-s}
\prod_{\fp|\ff,\ \fp\ndiv\tff}
(1-N\fp^{s-1}\sigma_{\fp,\tfm}\inv)\Phi_{\tfm}(s)
\eeq
For any abelian extension $K/k$ as above,
we now define a meromorphic, $\bbC G$-valued function by setting
\beql{eq:2C}
\Phi_{K/k}(s):=(|d_k|N\ff(K))^{s-1}\pi_{k(\fm(K)),K}(\Phi_{\fm(K)}(s))
\eeq
(This agrees with the case $T=\emptyset$ of equation~\cite[eq.~4]{sconzp}
although there $K$ was assumed to be totally real). Note in particular
the factor $(|d_k|N\ff(K))^{s-1}$ in this definition which,
when combined with\refeq{2B} gives the following analogue of\refeq{2A.75} in the same
situation
\beql{eq:piTheta}
\pi_{K,\tK}(\Phi_{K/k}(s))=
\prod_{\fp\in S_{\rm ram}(K/k)\atop \fp\nin S_{\rm ram}(\tK/k)}
(1-N\fp^{s-1}\sigma_{\fp,\tK}\inv)\Phi_{\tK/k}(s)
\eeq
Unlike $\Theta_{\fm}(s)$ and $\ThetaKk(s)$, however, there is no Euler
product for $\Phi_\fm(s)$ or $\PhiKks$ (but see
Thm.~\ref{thm:2A}).
The coefficient of $\sigma\inv$ in $\PhiKks$ may be denoted
$Z_{K/k}(s;\sigma)$ and loosely called `a twisted zeta-function of $K/k$'.

Now let $\chi:G\rightarrow\bbC^\times$ be any (irreducible) character of $G$
which may also be regarded as a character of $G_{\fm(K)}$ and hence
of $\Cl_{\fm(K)}(K)$. We write $\fm(\chi)=\ff(\chi)\fz(\chi)$ for the
conductor of the extension $K^{\ker(\chi)}/k$ cut out by
$\chi$. Then $\chi=\hat{\chi}\circ \pi_{\fm(K),\fm(\chi)}$ where
$\hat{\chi}$ is the \emph{primitive} character associated to $\chi$, as defined on
$\Cl_{\fm(\chi)}(k)$, hence on $G_{\fm(\chi)}$ and also on the group
of fractional ideals prime to $\ff(\chi)$. The corresponding $L$-function
is $L(s,\hat{\chi})$ which equals
$\prod_{\fp\ndiv\ff(\chi)}(1-N\fp^{-s}\hat{\chi}(\fp))\inv$ for ${\rm Re}(s)>1$.
Using this notation and extending $\chi$ linearly to $\bbC G$, equation\refeq{2A.5} gives
\beql{eq:2D}
\chi(\Theta_{K/k}(s))=\prod_{\fp|\ff(K),\ \fp\ndiv\ff(\chi)}
(1-N\fp^{-s}\hat{\chi}\inv(\fp)) L(s,\hat{\chi}\inv)
\eeq
An analogous equation for $\chi(\PhiKks)$
follows from~\cite[Thm.~3.3]{twizas} and\refeq{2C} above, namely
\beql{eq:2E}
\chi(\Phi_{K/k}(s))=g(\hat{\chi})(|d_k|N\ff(\chi))^{s-1}
\prod_{\fp|\ff(K),\ \fp\ndiv\ff(\chi)}
(1-N\fp^{s-1}\hat{\chi}\inv(\fp)) L(s,\hat{\chi})
\eeq
where $g(\hat{\chi})=g_{\fm(\chi)}(\hat{\chi})\in\barbbQ^\times$ is the Gauss
sum attached to $\hat{\chi}$ \emph{as a character of $\Cl_{\fm(\chi)}(k)$ or
$G_{\fm(\chi)}$}. (For a definition, see~\cite[\S 6.4]{twizas} or
our Remark~\ref{rem:2A}(ii)).
The last two equations give alternative definitions
of $\Theta_{K/k}(s)$ and $\Phi_{K/k}(s)$ respectively
since for any $x\in\bbC G$ we have $x=\sum_\chi\chi(x)e_\chi$
where $e_\chi:=|G|\inv\sum_{g\in G}\chi(g)g\inv$ is the idempotent of $\bbC G$
associated to $\chi$. However, for many purposes the
`equivariant' definitions\refeq{2Z} and\refeq{2A.9} are more helpful.

The function $\Theta_{K/k}$ -- and in particular, its r\^ole in the Stark
Conjectures -- is far better known than $\Phi_{K/k}$ so we now
give two relations between these functions.
For each place $v\in S_\infty$, we
write $c_v=c_{v,K}$ for
the unique generator of $D_v(K/k)$. Thus if $v|\fz(K)$  then $c_v$ is the associated
complex conjugation.
Otherwise (and in particular, if $v$ is itself complex) $c_v$ is trivial.
We also define a function $C_v:\bbC\rightarrow \bbC G$ by
\[
C_v(s)=\left\{
\begin{array}{ll}
e^{i\pi s}-e^{-i\pi s}c_v=2i\sin(\pi s)&
\mbox{if $v$ is complex}\\
e^{i\pi s/2}+e^{-i\pi s/2}c_v&
\mbox{if $v$ is real}
\end{array}
\right.
\]
\begin{thm}\label{thm:2A}
If $K$ is any abelian extension of $k$ then, with the above notations,
\[
i^{r_2(k)}\sqrt{|d_k|}\PhiKk(1-s)=
((2\pi)^{-s}\Gamma(s))^d\left(\prod_{v\in S_\infty}C_v(s)\right)\ThetaKk(s)
\]
as meromorphic, $\bbC G$-valued functions of $s\in\bbC$.
(Note that $i^{r_2(k)}\sqrt{|d_k|}$ is a square-root of $d_k$).
\end{thm}
\bPf\ It suffices to prove the $\chi$-part of this equation for
each character $\chi$. By means of equations\refeq{2D} and\refeq{2E} we are reduced to showing
\beql{eq:2F}
g(\hat{\chi})(|d_k|N\ff(\chi))^{-s}L(1-s,\hat{\chi})/L(s,\hat{\chi}\inv)=
((2\pi)^{-s}\Gamma(s))^d
i^{-r_2(k)}\sqrt{|d_k|}\inv
\prod_{v\in S_\infty}\chi(C_v(s))
\eeq
This is in fact the functional equation for the $L$-function. To put it into
a more familiar form we may use
the identities $\Gamma(z)\Gamma(1-z)=\pi/\sin(\pi z)$
and $\Gamma(z)\Gamma(z+\half)=(2\pi)^\half 2^{\half-2z}\Gamma(2z)$
to rewrite $\chi(C_v(s))$ as $2\pi i/(\Gamma(s)\Gamma(1-s))$ if $v$ is complex, as
$2^s\pi^\half\Gamma(s)\inv \left(i\Gamma(\frac{s+1}{2})/\Gamma(1-\frac{s}{2})\right)$
if $v$ is real and $v|\fz(\chi)$ (\ie\ $\chi(c_v)=-1$)
and as
$2^s\pi^\half\Gamma(s)\inv \left(\Gamma(\frac{s}{2})/\Gamma(1-\frac{s+1}{2})\right)$
if $v$ is real and $v\ndiv\fz(\chi)$ (\ie\ $\chi(c_v)=1$).
Equation\refeq{2F} then becomes
\begin{eqnarray*}
\lefteqn{
g(\hat{\chi})(|d_k|N\ff(\chi))^{-s}L(1-s,\hat{\chi})/L(s,\hat{\chi}\inv)=
2^{r_2(k)(1-2s)}\pi^{-d(s-\half)}\sqrt{|d_k|}\inv i^q\,\times
}\hspace{5em}&&\\
&&(\Gamma(s)/\Gamma(1-s))^{r_2(k)}(\Gamma((s+1)/2)/\Gamma((2-s)/2))^{q}
(\Gamma(s/2)/\Gamma((1-s)/2))^{r_1(k)-q}\\
\end{eqnarray*}
(where $q=|\fz(\chi)|=|\fz(\chi\inv)|$) and this, after rearranging,
is the functional equation in the form
of~\cite[Eq. (5)]{Tatu}, with `$\chi$' and `$\bar{\chi}$' replaced by our
$\hat{\chi}\inv$ and $\hat{\chi}$ respectively. (The Gauss sum which would
be denoted $F(\hat{\chi}\inv)$ in~\cite{Tatu}
is our $g(\hat{\chi})$ by~\cite[Rem. 6.3]{twizas}.)\ePf
\rem\ One could of course use Theorem~\ref{thm:2A} to \emph{define} $\PhiKk(s)$
in terms of $\ThetaKk(s)$ but the definition given via
$\Phi_\fm(s)$ and the twisted zeta-functions seems a little more natural. In any case,
$\Phi_\fm(s)$ has several important properties which we shall exploit in our
study of $\PhiKk(s)$, especially in Theorem~\ref{thm:3A}.\vspace{1ex}\\
For any abelian group $H$ and commutative ring $R$
we define an involutive automorphism $\underline{\ }^\ast$ of $RH$ by setting
$(\sum a_hh)^\ast=\sum a_hh\inv$.
Our second relation expresses $\Phi_{\fm}(s)^\ast$ in terms of the functions
$\Theta_{\fg\fz}(s)$ as
$\fg$ runs over the ideals dividing $\ff$.
For the rest of this section
we shall use the notation $\fn$ generically to represent the cycle $\fg\fz$
for an ideal $\fg$ dividing $\ff$.
Let $\tilde{\nu}_{\fn,\fm}:\bbC G_{\fn}\rightarrow \bbC G_\fm$ be
the linear extension of the map sending an element of $G_{\fn}$ to the sum of its
pre-images under $\pi_{\fm,\fn}$. Thus $\nu_{\fn,\fm}:=|\ker(\pi_{\fm,\fn})|\inv
\tilde{\nu}_{\fn,\fm}$ is a ring homomorphism right inverse to $\pi_{\fm,\fn}$.
As in~\cite{twizas}, we write
$k^\times_{\fn}$ for the group of elements
of $k$ which are congruent to $1$~$\mbox{mod}^\times\,\fn$ and $E_{\fn}$ for
$E(k)\cap k_{\fn}^\times$. For any fractional ideal $I$ we denote by
$\cT(\fg,I)$ the set of precise $\fg$-torsion classes in the $\cO$-module
$k/I$ (in particular $\cT(\fg,I)\subset {\fg}\inv I/I$).
These sets parametrize $\Cl_{\fn}(k)$ as explained
in~\cite[Prop./Def. 6.1]{twizas}. In particular, each $y\in\cT(\fg,\fg J\inv)$
gives a well-defined class $[y;J]_{\fn}\in\Cl_{\fn}(k)$
lying in the fibre  of $\pi_{\fn,\fz}$ over
$[J]_\fz$ (which we  denote
$\pi_{\fn,\fz}\inv([J]_\fz)$)
and defined by $[y;J]_{\fn}:=[bJ]_{\fn}$ for any $b$ in $y\cap k^{\times}_\fz$.
We write $\be(z)$ for $\exp(2\pi i z)$, $\Tr$ for $\Tr_{k/\bbQ}$ and
$\cD_k$ for the different of $k/\bbQ$. Thus the map
$u\mapsto \be(\Tr(u))$ is a well-defined (additive) character of $k/\cD_k\inv$, taking
${\fg}\inv\cD_k\inv/\cD_k\inv$ onto $\mu_{f'}$ where  $f'\in\bbZ$ denotes the
positive generator of the ideal $\fg\cap \bbZ$ of $\bbZ$. Let
\[
A_{\fn}:=\sum_{u\in\cT(\fg,\cD_k\inv)}
\be(\Tr(u))
\sigma_{[u;\fg\cD_k]_{\fn}}\in\bbZ[\mu_{f'}]G_{\fn}
\]
\begin{thm}\label{thm:2B}
With notations and hypotheses as above,
\begin{eqnarray}
\Phi_{\fm}(s)^\ast&=&
  \sum_{\fg|\ff}\left(\frac{N\ff}{N\fg}\right)^{-s}[E_{\fn}:E_{\fm}]
  \tilde{\nu}_{\fn,\fm}(A_{\fn}\Theta_{\fn}(s))\label{eq:2F.5}\\
                  &=&
  \sum_{\fg|\ff}\left(\frac{N\ff}{N\fg}\right)^{1-s}
  \left(\prod_{\fp|\ff,\ \fp\ndiv\fg}(1-N\fp\inv)\right)
  \nu_{\fn,\fm}(A_{\fn}\Theta_{\fn}(s))\label{eq:2F.75}
\end{eqnarray}
\end{thm}
\rem\label{rem:2A}

(i) If $\fz=\emptyset$ then $[-u;\fg\cD_k]_{\fn}=[u;\fg\cD_k]_{\fn}$.
It follows that the coefficients of $A_{\fn}$ are real in this case.

(ii) For any character $\chi:G_\fm\rightarrow \bbC^\times$ the Gauss sum
$g(\chi)=g_{\fm}(\chi)$ equals $\chi\inv(A_\fm)$.

(iii) A rather complicated `functional equation for $\Theta$' might be obtained
by combining Theorem~\ref{thm:2B} (with $1-s$ for $s$) and Theorem~\ref{thm:2A}.
\vertsp\\
\noindent \textsc{Proof of Thm.\ 1.2\ }\
By meromorphic continuation we can assume ${\rm Re}(s)>1$
so that all the sums below are absolutely convergent. The definition of
the twisted zeta-function $Z_{\emptyset}$ in~\cite{twizas} together with
Lemma~6.1 of~\emph{ibid.} with $H=E_\fm$, $T=\emptyset$ yields the expresssion
\[
\Phim(s)^\ast=[E_\fz:E_\fm]\inv
\sum_{J\in\cJ}
\sum_{y\in\cT(\ff,\ff J\inv)}
\left(
\sum_{a\in(\ff\inv\cD_k\inv J\cap k_\fz^\times)/E_\fm}
\be(\Tr(ay))N(a\ff\cD_k J\inv)^{-s}
\right)
\sigma_{[y;J]_\fm}
\]
where $\cJ$ is any set of fractional ideals representing
$\Cl_\fz(k)$ and $a$ runs over a set of representatives for
the action of $E_\fm$ on $\ff\inv\cD_k\inv J\cap k_\fz^\times$. Now, for fixed
$y$, the value of $\be(\Tr(ay))$ depends only on the class $w$
of $a$ in the $\cO$-module $\ff\inv\cD_k\inv J/\cD_k\inv J$. Writing $\ff\inv\cD_k\inv J$ as
the union of such classes $w$ and grouping them according to
their $\cO$-annihilator $\fg$ (which must divide $\ff$) we see
that the sum over $a$ in the last equation may be written as
\begin{eqnarray*}
\lefteqn{\sum_{\fg|\ff}N(\ff{\fg}\inv)^{-s}\sum_{w\in \cT(\fg,\cD_k\inv J)}\be(\Tr(wy))
\sum_{(a\in w\cap k_\fz^\times)/E_\fm}N(a\fg\cD_k J\inv)^{-s}=}\hspace*{10ex}&&\\
&&\sum_{\fg|\ff}N(\ff{\fg}\inv)^{-s}\sum_{w\in \cT(\fg,\cD_k\inv J)}\be(\Tr(wy))
[E_{\fn}:E_\fm]\zeta_{\fn}(s,[w;\fg\cD_k J\inv]_{\fn})
\end{eqnarray*}
by part~(iv) of~\cite[Prop./Def.~6.1]{twizas}
applied to $w$ in place of $y$, \emph{mutatis mutandis}. Hence
\begin{eqnarray}
\lefteqn{\Phim(s)^\ast=
\sum_{\fg|\ff}N(\ff{\fg}\inv)^{-s}[E_\fz:E_{\fn}]\inv\times}\hspace*{5ex}&&\nonumber\\
&&\sum_{J\in\cJ}
\left[
\sum_{w\in\cT(\fg,\cD_k\inv J)}
\zeta_{\fn}(s,[w;\fg\cD_k J\inv]_{\fn})
\sum_{y\in\cT(\ff,\ff J\inv)}\be(\Tr(wy))\sigma_{[y;J]_\fm}
\right]\label{eq:2G}
\end{eqnarray}
Now for fixed $\fg$ and $J$,
the reduction map $\delta:J\inv/\ff J\inv\rightarrow J\inv/\fg J\inv$
takes $\cT(\ff,\ff J\inv)$ onto $\cT(\fg,\fg J\inv)$.
Let us gather up the terms in the last sum of\refeq{2G} according to the
value of $z=\delta(y)$. We note that
$\be(\Tr(wy))=\be(\Tr(wz))$ and we claim that
for each $z\in\cT(\fg,\fg J\inv)$,
\beql{eq:2H}
\sum_{{y\in\cT(\ff,\ff J\inv)}\atop{\delta(y)=z}}\sigma_{[y;J]_\fm}=
[E_{\fn}:E_{\fm}]\tilde{\nu}_{\fn,\fm}(\sigma_{[z;J]_{\fn}})
\eeq
(proof deferred). Furthermore,
if $z$ is fixed, then as $w$ ranges through $\cT(\fg,\cD_k\inv J)$
so $u:=wz$ ranges exactly once through $\cT(\fg,\cD_k\inv)$ and
$[u;\fg\cD_k]_{\fn}=[w;\fg\cD_k J\inv]_{\fn}[z;J]_{\fn}$
(see Props.~6.1 and~6.2~(i) of~\cite{twizas}).
Therefore the term in square brackets in\refeq{2G} can be rewritten as
\begin{eqnarray}
\lefteqn{\sum_{w\in\cT(\fg,\cD_k\inv J)}
\zeta_{\fn}(s,[w;\fg\cD_k J\inv]_{\fn})
\sum_{z\in\cT(\fg,\fg J\inv)}
[E_{\fn}:E_{\fm}]\be(\Tr(wz))\tilde{\nu}_{\fn,\fm}(\sigma_{[z,J]_{\fn}})
=
}\hspace*{3ex}&&\nonumber\\
&&[E_{\fn}:E_{\fm}]
\tilde{\nu}_{\fn,\fm}\left(
\sum_{z\in\cT(\fg,\fg J\inv)}
\left\{\sum_{u\in\cT(\fg,\cD_k\inv)}
\be(\Tr(u))\zeta_{\fn}(s,[u;\fg\cD_k]_{\fn}[z;J]_{\fn}\inv)
\right\}
\sigma_{[z,J]_{\fn}}
\right)\nonumber
\end{eqnarray}
Notice that the term in braces is simply the coefficient of
$\sigma_{[z,J]_{\fn}}$ in $A_{\fn}\Theta_{\fn}(s)$. Note also
as $J$ runs through $\cJ$ and $z$ through $\cT(\fg,\fg J\inv)$ for each $J$,
so $\sigma_{[z,J]_{\fn}}$ runs exactly $[E_\fz:E_{\fn}]$ times through
$G_{\fn}$ (by~\cite[Prop./Def.~6.1]{twizas} again).
Therefore, substituting the R.H.S.
of the last equation for the term in square brackets in\refeq{2G}, the factor
$[E_\fz:E_{\fn}]\inv$ is cancelled, giving\refeq{2F.5}.
Equation\refeq{2F.75} follows easily from the latter
and the equality $[E_{\fn}:E_{\fm}]|\ker(\pi_{\fm,\fn})|=
|(\cO/\ff)^\times|/|(\cO/\fg)^\times|$ which
follows  in turn from the exact sequence
\beql{eq:2I}
1\rightarrow E_\fz/E_\fm\longrightarrow
(\cO/\ff)^\times\longrightarrow\Cl_{\fm}(k)
\longrightarrow\Cl_\fz(k)\rightarrow 1
\eeq
and a similar one with $\fn$ and $\fg$ in place of $\fm$ and
$\ff$.
It only remains to establish equation\refeq{2H}. But this follows from the commutativity
of the square\refeq{2J} of surjective maps
\begin{figure}
\setlength{\unitlength}{4mm}
\beql{eq:2J}
\begin{picture}(10,8)(-1,0)
\put(2.5,0.5){\vector(1,0){3}}
\put(2.5,7.5){\vector(1,0){3}}
\put(0,6.7){\vector(0,-1){5.4}}
\put(8,6.7){\vector(0,-1){5.4}}
\put(4,1){\makebox(0,0)[b]{$\beta$}}
\put(4,8){\makebox(0,0)[b]{$\alpha$}}
\put(0.5,4){\makebox(0,0)[l]{$\delta$}}
\put(8.5,4){\makebox(0,0)[l]{$\pi_{\fm,\fn}$}}
\put(0,0.5){\makebox(0,0){$\cT(\fg,\fg J\inv)$}}
\put(0,7.5){\makebox(0,0){$\cT(\ff,\ff J\inv)$}}
\put(8,0.5){\makebox(0,0){$\pi_{\fn,\fz}\inv([J]_\fz)$}}
\put(8,7.5){\makebox(0,0){$\pi_{\fm,\fz}\inv([J]_\fz)$}}
\end{picture}
\eeq
\end{figure}
(with $\alpha(y):=[y;J]_\fm$, $\beta(z)=:[z;J]_{\fn}$) together with the following lemma.
\begin{lemma}
Given $z\in \cT(\fg,\fg J\inv)$ and $\fc\in \pi_{\fm,\fz}\inv([J]_\fz)$ such that
$\beta(z)=\pi_{\fm,\fn}(\fc)$, there exists precisely one $E_{\fn}$-orbit of elements
$y\in \cT(\ff,\ff J\inv)$ such that $\delta(y)=z$ and $\alpha(y)=\fc$. In particular,
there are exactly $[E_{\fn}:E_\fm]$ such elements.
\end{lemma}
This lemma is in turn easily deduced from~\cite[Prop./Def.~6.1]{twizas}
with some diagram chasing.\ePf
\noindent Now suppose that $K/k$ is an abelian extension with group $G$ such that $\fm(K)=\fm$.
Then we can use Theorem~\ref{thm:2B} to obtain rather complicated expressions
for $\PhiKk(s)^\ast$ in terms of the $\Theta_{\tK/k}(s)$ as $\tK/k$ runs over certain sub-extensions
of $K/k$. As an example we record the specialisation at $s=0$ of one
such expression which will be useful in the next section and, potentially, for
computation.
For each $\fn|\fm$ as in Theorem~\ref{thm:2B},
we set $K[\fn]=K\cap k(\fn)$ and for any $\tilde{K}$ with $K\supset \tilde{K}\supset k$
we define $\tilde{\nu}_{\tilde{K},K}:\bbC \Gal(\tilde{K}/k)\rightarrow \bbC G$ in a manner
entirely analogous to $\tilde{\nu}_{\fn,\fm}$. It is easy to see that
\[
\pi_{k(\fm),K}\circ\tilde{\nu}_{\fn,\fm}=[k(\fm):Kk(\fn)]\tilde{\nu}_{K[\fn],K}
\circ \pi_{k(\fn),K[\fn]}
\]
as maps from $\bbC G_{\fn}$ to $\bbC G$.
Therefore, equations\refeq{2F.5} and\refeq{2C} give
\[
\PhiKko^\ast=\frac{1}{|d_k|N\ff}\sum_{\fg|\ff}[E_{\fn}:E_{\fm}][k(\fm):Kk(\fn)]
\tilde{\nu}_{K[\fn],K}
(\pi_{k(\fn),K[\fn]}(A_{\fn}\Theta_{\fn}(0)))
\]
Since $\fm(K[\fn])$ divides -- but is not in general equal to -- $\fn$, we use\refeq{2A}
and\refeq{2A.25} to calculate $\pi_{k(\fn),K[\fn]}(\Theta_{\fn}(0))$ giving
\begin{cor}\label{cor:2A} Suppose $\fm(K)=\fm$ then with the above notation
($\fm=\fg\fz$ \etc) we have
\beql{eq:2J.5}
\PhiKko^\ast=\frac{1}{|d_k|N\ff}\sum_{\fg|\ff}
\tilde{\nu}_{K[\fn],K}
(B_{\fn}\Theta_{K[\fn]/k}(0))
\eeq
where, for each $\fg|\ff$, the element $B_{\fn}$ of $\bbZ[\mu_{f'}]\Gal(K[\fn]/k)$
is given by
\[
B_{\fn}=[E_{\fn}:E_{\fm}][k(\fm):Kk(\fn)]\pi_{k(\fn),K[\fn]}(A_{\fn})
\prod_{\fp|\fn\atop\fp\ndiv\fm(K[\fn])}(1-\sigma_{\fp,K[\fn]}\inv)
\]
\ePf
\end{cor}
\section{The Behaviour of $\Phi_{K/k}(0)$}\label{sec:behav}
The main object of interest in this paper is the value of $\Phi_{K/k}(s)$ at $s=0$,
particularly from a $p$-adic viewpoint.
We shall suppose to start with that $K/k$ is any abelian extension
with notation as in the previous section.
For any character $\chi:G\rightarrow\bbC^\times$,
equation\refeq{2E} shows that the order of
vanishing of $\chi(\PhiKks)$ at $s=0$ is the same as that of
$L(s,\hat{\chi})$. One knows (\eg\ by the functional equation) that for
$\chi\neq\chi_0$ (the trivial character of $G$) the latter order
equals the number of places $v\in S_\infty$  for which
$D_v(K/k)\subset\ker(\chi)$ while if $\chi=\chi_0$ then it equals
$|S_\infty|-1$. A first consequence is that $\PhiKko$ vanishes
unless $k$ is either a totally real \emph{or} an imaginary quadratic
field. Moreover, in the latter case it is given as an explicit rational multiple of
$\zeta_k(0)e_{\chi_0}=-(h_k/w_k)e_{\chi_0}$ by\refeq{2E}
and is of no great interest in the present context.
\emph{We shall therefore assume henceforth the following}
\begin{hyp}
The base field $k$ is totally real.
\end{hyp}
This condition implies that $d_k$ is a positive integer.
Siegel-Klingen and Shintani's results mentioned in the introduction,
imply that $\Theta_{L/k}(m)$ has rational coefficients for any abelian
$L/k$ and any $m\in\bbZ_{\leq 0}$. It follows from Cor.~\ref{cor:2A}, that $\PhiKko$ lies in
$\bbQ(\mu_{f(K)})G$ where $f(K)$ is the positive generator of $\ff(K)\cap
\bbZ$ (see also~\cite[Lemma 3.3]{zetap1}).
Let $e^-_{K/k}$ denote the idempotent $\prod_{v\in S_\infty}(\half(1-c_v))$ of $\bbQ G$.
\begin{prop}\label{prop:2A} If $k\neq \bbQ$
then $\PhiKko$ is a generator of the ideal $e^-_{K/k}\bbQ(\mu_{f(K)})G$
of $\bbQ(\mu_{f(K)})G$. For $k=\bbQ$ the same is true of
$\Phi_{K/\bbQ}(0)+\half\prod_{q|f(K)}(1-q\inv)e_{\chi_0}$.
\end{prop}
\bPf\ Except in the case $k=\bbQ$ and
$\chi=\chi_0$, the previous discussion gives the equivalences
$\chi(\PhiKko)\neq 0\Leftrightarrow\chi(c_v)\neq 1\ \forall v
\Leftrightarrow\chi(c_v)=-1\ \forall v$.
The Proposition follows easily, using equation\refeq{2E}
to show that $\chi_0(\Phi_{K/\bbQ}(0))=-\half\prod_{q|f(K)}(1-q\inv)$.
\ePf
\rem\label{rem:2A.5}

(i) Let $H^+$  and $H^-$ be the subgroups of $G$ generated by the sets
$\{c_v:v\in S_\infty\}$ and $\{c_vc_{v'}:v,v'\in S_\infty\}$ respectively.
Thus $G\supset H^+\supset H^-$ and the index $|H^+:H^-|$ is $1$ or $2$.
Now $c_ve^-_{K/k}=-e^-_{K/k}\ \forall\,v\in S_\infty$ and it follows that
$H^-$ fixes $e^-_{K/k}$. Thus, if $H^+=H^-$ ($\Leftrightarrow 1$ is the product of an
odd number of $c_v$'s) then $e^-_{K/k}$
vanishes and so will $\PhiKk(0)$, by Proposition~\ref{prop:2A}, unless $k=\bbQ$.
Thus one loses very little in assuming that $|H^+:H^-|=2$. This condition is equivalent
to the statement that $K$ contains a $CM$ subfield and so, in particular, is totally complex.
Indeed the unique maximal CM subfield is $K^-:=K^{H^-}$ and $K^+:=K^{H^+}$ is its
maximal real subfield. In this case, one can show that $e^-_{K/k}$ is non-zero.
The Proposition therefore implies that $\PhiKko$ is non-zero and fixed by $H^-$.
It follows from\refeq{piTheta} (with $\tilde{K}=K^-$) that in this case
\beql{eq:hey!}
\PhiKko=
|H^-|\inv\tilde{\nu}_{K^-,K}(
\prod_{\fp\in S_{\rm ram}(K/k)\atop \fp\nin S_{\rm ram}(K^-/k)}
(1-N\fp^{-1}\sigma_{\fp,K^-}\inv)
\Phi_{K^-/k}(0)
)
\eeq
Furthermore, $|H^-|$ divides $2^{d-1}$. For certain purposes, this
formula allows one to further reduce to the case where $K$ is a CM field.

(ii) Using equation\refeq{2D} in place of\refeq{2E} one can show similarly that
$\ThetaKk(0)=0$ unless $k$ is totally real or $k$ is imaginary quadratic and
$K/k$ is unramified. Assuming the former, one finds that $\ThetaKk(0)$ lies in
$e^-_{K/k}\bbQ G$ (unless $K=k=\bbQ$) and the analogue of\refeq{hey!} holds with
$\Theta$ in place of $\Phi$ and the product replaced by
$\prod(1-\sigma_{\fp,K^-}\inv)$. Since the latter lies in  $\bbZ G$, the reduction
to the CM case is even easier than for $\PhiKk(0)$ (\cf~\cite[Thm.~IV.5.2]{Tate}).

(iii) In general, $\Theta_{K/k}(s)$ `has more zeroes' at $s=0$ than
$\Phi_{K/k}(s)$ \ie\ $\ThetaKk(0)$ may not generate the whole of
$e^-_{K/k}\bbQ G$ over $\bbQ G$.
More precisely, a comparison of equations\refeq{2D} and\refeq{2E}
shows that $\chi(\Theta_{K/k}(0))$ vanishes
not only when $\chi(\PhiKko)$ does but also when
$D_{\fp}(K/k)\subset\ker(\chi)$ for some \emph{finite} prime
$\fp\in S_{\rm ram}(K/k)$.
This also explains the necessity of using not only $\Theta_{K/k}(0)$ but also
$\Theta_{\tK/k}(0)$ (for \emph{sub}fields $\tK$ of $K$) to express $\PhiKko$, as
in \eg\ Cor.~\ref{cor:2A}.

(iv) On the other hand, Theorem~\ref{thm:2A} gives a simple relation between $\PhiKko$
and $\Theta_{K/k}(s)$ at $s=1$. To avoid some relatively unimportant complications coming
from the pole of the latter, we write $\Theta^{\rm n.t.}_{K/k}(s)$ for
the function $(1-e_{\chi_0})\ThetaKk(s)$ which is regular at $s=1$. Then
Theorem~\ref{thm:2A}  gives
\[
\sqrt{d_k}(1-e_{\chi_0})\PhiKko=(i/\pi)^d e^-_{K/k}\Theta^{\rm n.t.}_{K/k}(1)
\]
In particular, the coefficients of $e^-_{K/k}\Theta^{\rm n.t.}_{K/k}(1)$ are algebraic numbers
multiplied by $\pi^d$
By contrast, the coefficients of $(1-e^-_{K/k})\Theta^{\rm n.t.}_{K/k}(1)$
are expected to contain other transcendental terms (as well as powers of $\pi$)
which encode arithmetic information on $K$.
Indeed, if $\chi$ is neither totally odd nor trivial then
the Stark Conjectures predict factors in
$\chi(\Theta^{\rm n.t.}_{K/k}(1))$
coming from logarithms of absolute values of
units of $K$. (See~\cite{Tate}, particularly Conjecture~I.8.2
and Lemme~I.8.7).
\vspace{1ex}\\
Every element $\alpha$ of $\Gal(\barbbQ/\bbQ)$
induces an automorphism of $\bbQ(\mu_{f(K)})G$ by its action on coefficients.
To describe its effect on $\PhiKko$, we
let $\bbQ^{\rm ab}$ and $k^{\rm ab}$ denote the maximal abelian extensions of
$\bbQ$ and $k$ respectively inside $\barbbQ$ and we write
${\rm Ver}$ for the transfer homomorphism from $\Gal(\bbQ^{\rm ab}/\bbQ)$
(identified with the abelianisation $\Gal(\barbbQ/\bbQ)^{\rm ab}$) to
$\Gal(k^{\rm ab}/k)$ (identified with $\Gal(\barbbQ/k)^{\rm ab}$).
If $F$ is any extension of $k$ within
$k^{\rm ab}$ we compose ${\rm Ver}$ with the restriction
to get a homomorphism $\cV_F:\Gal(\bbQ^{\rm ab}/\bbQ)\rightarrow\Gal(F/k)$.
For $\alpha$ as above, it follows from~\cite[Prop.~3.1]{zetap1}
that $\alpha(\Phi_{\fm(K)}(0))=\cV_{k(\fm(K))}(\alpha|_{\bbQ^{\rm ab}})\Phi_{\fm(K)}(0)$
(product in $\bbQ(\mu_{f(K)})G_{\fm(K)}$ on the R.H.S.). Applying
$(d_k N\ff(K))\inv\pi_{k(\fm(K)),K}$ to both sides gives
\begin{prop}\label{prop:2B} For each $\alpha\in\Gal(\barbbQ/\bbQ)$ we have
$\alpha(\PhiKko)=\cV_{K}(\alpha|_{\bbQ^{\rm ab}})\PhiKko$ (product in $\bbQ(\mu_{f(K)})G$).\ePf
\end{prop}
We now turn to the $p$-adic properties of $\PhiKko$ where
$p$ is a prime number, fixed throughout this paper. We fix also an embedding
$j:\bar{\bbQ}\rightarrow\bar{\bbQ}_p$ which extends
coefficientwise to a homomorphism of group rings
$j:\bar{\bbQ}H\rightarrow\bar{\bbQ}_pH$ for any group $H$.
In particular, for any abelian extension $K$ of $k$,
$j(\PhiKko)$ is an element of $j(\bbQ(\mu_{f(K)}))G\subset\bbQ_p(\mu_{f(K)})G$.

Let us define the $p$-adic regulator $R_{K/k,p}=R^{(j)}_{K/k,p}$ mentioned in the
introduction. We use the natural topological ring isomorphism to identify
$\bbQ_p\otimes_\bbQ K$ with the product
$\prod_{\fP|p}K_\fP$ of the completions of $K$ at the primes $\fP$ of $\cO_K$ dividing $p$.
This allows us to regard each $K_\fP^\times$ as a \emph{sub}group
of $(\bbQ_p\otimes_\bbQ K)^\times$ and identifies
$(\bbZ_p\otimes_\bbZ\cO_K)^\times$ (\emph{resp.} its Sylow pro-$p$ subgroup) with
the $p$-semilocal unit group $U_p(K):=\prod_{\fP|p}U(K_\fP)$
(\emph{resp.} with $U^1_p(K):=\prod_{\fP|p}U^{1}(K_\fP)$).
Each $g\in G$ acts by $1\otimes g$ on
$(\bbQ_p\otimes_\bbQ K)^\times$,
sending the subgroup $K_\fP^\times$ isomorphically onto $K_{g(\fP)}^\times$ for each
$\fP|p$ and making $U^1_p(K)$ into a natural $\bbZp G$-module.

For each $i=1,\ldots,n$, we write $\fP_i$ for the prime of $K$ above $p$
corresponding to the embedding $j\tau_i:=j\circ\tau_i:K\rightarrow\barbbQp$.
namely $\fP_i=\{a\in\cO_K:|j\tau_i(a)|_p<1\}$and we write $\fp_i$ for $\fP_i\cap\cO$.
(Note that $\{\fp_1,\ldots,\fp_d\}=\{\fp|p\}$ but in general we may have $\fp_i=\fp_{i'}$, and even $\fP_i=\fP_{i'}$,
for $i\neq i'$.) The the embedding $j\tau_i$
extends to a homomorphism $\bbQ_p\otimes_\bbQ K\rightarrow\barbbQp$
which is the composite of the projection onto $K_{\fP_i}$
with the natural isomorphism $K_{\fP_i}\rightarrow\overline{j(\tau_i(K))}$
(topological closure in $\bbC_p$) induced by $j\tau_i$.
Both this extension of $j\tau_i$ and this isomorphism will, by abuse,
also be denoted $j\tau_i$.
We define a $\bbZ G$-homomorphism
$\lambda_{i,p}=\lambda^{(j)}_{K/k,i,p}:(\bbQ_p\otimes_\bbQ K)^\times\rightarrow\bbC_p G$ by setting
$\lambda_{i,p}(a)=\sum_{g\in G}\log_p(j\tau_i(ga))g\inv$ where
$\log_p$ denotes Iwasawa's $p$-adic logarithm.
Thus $\lambda_{i,p}$ factors through the projection onto $\prod_{\fP|\fp_i}K_\fP^\times$
and actually takes values in $\overline{j(\tau_i(K))}G$. Furthermore, its restriction to
$U_p(K)$ also factors through the natural projection of this group onto $U^1_p(K)$ and
its restriction to the latter is $\bbZ_p G$-linear. There is therefore a unique
$\bbZ G$-homomorphism $R_{K/k,p}^{(j)}:\bigwedge_{\bbZ G}^dU_p(K)\rightarrow\barbbQ_p G$
taking $u_1\wedge\ldots\wedge u_d$ to $\det(\lambda_{i,p}(u_t))_{i,t=1}^d$.
Moreover $R^{(j)}_{K/k,p}$ is the composite of the natural map
$\bigwedge_{\bbZ G}^dU_p(K)\rightarrow\bigwedge_{\bbZ_p G}^dU^1_p(K)$
with a $\bbZ_p G$-homomorphism
$\bigwedge_{\bbZ_p G}^dU^1_p(K)\rightarrow\barbbQ_p G$ which is also denoted $R^{(j)}_{K/k,p}$
and is determined by the same formula.\vspace{1ex}\\
\rem\label{rem:imR}
The precise kernel and image of $R^{(j)}_{K/k,p}$ are hard to determine but
it is not difficult to see that the former is finite
and the latter spans $\barbbQ_pG$ over $\barbbQ_p$.
We sketch a proof. For the image, it obviously
suffices to construct $\theta_0\in\bigwedge_{\bbZ_p G}^dU^1_p(K)$ such
that $R^{(j)}_{K/k,p}(\theta_0)$ lies in $\barbbQ_pG^\times$.
Now, by means of the the $p$-adic exponential function in each completion $K_\fP$,
one can construct ${\rm Exp}_p:p\cO_K\rightarrow U^1_p(K)$ such that
$\lambda_{i,p}({\rm Exp}_p(x))=\rho_i(x):=\sum_{g\in G}j\tau_i(gx)g\inv$
(a sort of resolvent of $x$). Let $\theta_0=u_1\wedge\ldots\wedge u_d$
where $u_t={\rm Exp}_p(y_t\alpha)$, $\{y_1,\ldots,y_d\}$ is any $\bbQ$-basis
of $k$ and $\alpha$ generates a normal basis
for $K$ over $k$ (\ie\ $K=kG\alpha$),
these being chosen so that $y_t\alpha\in p\cO_K\ \forall\,t$.
Then $R^{(j)}_{K/k,p}(\theta_0)=\det(\rho_i(y_t\alpha))_{i,t=1}^d$.
That this lies in $\barbbQ_p G^\times$ follows from a simple calculation, the fact that
$\sum_{g\in G}g(\alpha) g\inv$ lies in $KG^\times$, see~\cite[Prop.~I.4.1]{Fro} \etc
For the kernel, one shows easily that ${\rm Exp}_p$ extends to a $\bbZ_p G$-linear
injection on $p\cO_K\otimes_\bbZ\bbZ_p$ and then that
the $u_t$ for $t=1,\ldots,d$ span a free $\bbZ_pG$-submodule of rank $d$ whose index
in $U^1_p(K)$ is \emph{finite}. The finitude of $\ker R^{(j)}_{K/k,p}$ now follows from the
fact that $R^{(j)}_{K/k,p}(\theta_0)$ is a unit.\vspace{1ex}\\
The dependence of $R^{(j)}_{K/k,p}$ on $j$
is determined as follows.
\begin{prop}\label{prop:2C}
For any $\alpha\in\Gal(\barbbQ/\bbQ)$ and $\theta\in \bigwedge_{\bbZ G}^dU_p(K)$,
we have
\[
(j\circ\alpha)(\sqrt{d_k})R^{(j\circ\alpha)}_{K/k,p}(\theta)=
\cV_{K}(\alpha|_{\bbQ^{\rm ab}})j(\sqrt{d_k})R^{(j)}_{K/k,p}(\theta)
\]
Moreover $j(\sqrt{d_k})R^{(j)}_{K/k,p}(\theta)$ lies in $\bbQ_p(\mu_{f(K)})G$.
\end{prop}
\bPf\ The proof of the first statement uses the definition of the
transfer map and differs in little but notation from that of the
second statement of~\cite[Prop.~3.2]{zetap1}, to which we refer.
The second statement above follows from the first,
using also the class field theoretic fact that $\cV_{K}=\pi_{k(\fm(K)),K}
\circ \cV_{k(\fm(K))}$ factors through the restriction map
$\Gal(\bbQ^{\rm ab}/\bbQ)\rightarrow\Gal(\bbQ(\mu_{f(K)})/\bbQ)$
(see \emph{loc. cit.} for details, noting that there $K$ is assumed to be a
totally real ray-class field).\ePf
\begin{prop}\label{prop:2D} For any $\theta\in \bigwedge_{\bbZ G}^dU_p(K)$ the element
$j(\sqrt{d_k}\PhiKko^\ast)R^{(j)}_{K/k,p}(\theta)$ of $\bar{\bbQ}_p G$
actually has coefficients in $\bbQ_p$
and is independent of the choice of the embedding $j$.
\end{prop}
\bPf\ Propositions~\ref{prop:2B} and~\ref{prop:2C} show that
$j(\sqrt{d_k}\PhiKkoa) R^{(j)}_{K/k,p}(\theta)$ is unchanged when $j$ is replaced by
$j\circ\alpha$ for any $\alpha\in\Gal(\barbbQ/\bbQ)$ which proves the second statement.
For the first, note that any $\beta\in\Gal(\barbbQ_p/\bbQ_p)$
commutes with $\log_p$. Hence, letting it act on coefficients we find:
\[
\beta(j(\sqrt{d_k}\PhiKkoa) R^{(j)}_{K/k,p}(\theta))=
(\beta\circ j)(\sqrt{d_k}\PhiKkoa) R^{(\beta\circ j)}_{K/k,p}(\theta)=
j(\sqrt{d_k}\PhiKkoa) R^{(j)}_{K/k,p}(\theta)
\]
by the second statement. Now let $\beta$ run through $\Gal(\barbbQ_p/\bbQ_p)$.
\ePf
\noindent\ As a consequence, the following makes sense and is independent of
the choice of $j$ above:
\begin{defn}\label{def:sKkSKk} For $K/k$ as above, we define
$\fs_{K/k}=\fs_{K/k,p}:\bigwedge_{\bbZ G}^dU_p(K)\rightarrow\bbQ_pG$ by setting
$\fs_{K/k}(\theta):=j(\sqrt{d_k}\PhiKkoa) R^{(j)}_{K/k,p}(\theta)$
for every $\theta\in\bigwedge_{\bbZ G}^dU_p(K)$.
We write $\fS_{K/k}=\fS_{K/k,p}$ for the image of $\fs_{K/k}$.
\end{defn}
Thus $\fs_{K/k}$ is $\bbZ G$-linear and factors through
a $\bbZ_p G$-linear map on $\bigwedge_{\bbZ_p G}^dU^1_p(K)$.
Consequently its image $\SKk$ is a (finitely generated)
$\bbZ_p G$-submodule of $\bbQ_pG$.\vspace{1ex}\\
\rem\label{rem:another}

(i) Slight variants of the proofs of Props.~\ref{prop:2C} and~\ref{prop:2D}
give the following: If $x_t\in K$ for $t=1,\ldots,d$ then
$\sqrt{d_k}\PhiKko^\ast\det(\sum_{g\in G}\tau_ig(x_t)g\inv)_{i,t}$ lies in $\bbQ G$.
This result illuminates Prop.~\ref{prop:2D} which may be deduced from it as follows:
take $x_t$ to be the $N$th truncation of the logarithmic series
evaluated at an element $u_t$ of $K \cap U^1_p(K)$ (which is dense in $U^1_p(K)$).
Then apply $j$ to the result and let $N\rightarrow\infty$.

(ii) Proposition~\ref{prop:2A} and Remark~\ref{rem:imR} imply that the $\bbQ_p$-span of
$\SKk$ is equal to $e_{K/k}^-\bbQ_p G$ or to $e_{K/k}^-\bbQ_p G+\bbQ_pe_{\chi_0}$
according as $k\neq\bbQ$ or $k=\bbQ$. Similarly they determine the kernel of
$\sKk$ up to finite index.

(iii) We note that the maps $\lambda^{(j)}_{K/k,i,p}$, hence
also $R^{(j)}_{K/k,p}$ and $\fs_{K/k}$, depend
on the choice and/or ordering of the coset representatives
$\tau_i,\ldots,\tau_d$. However, changing these choices clearly only
multiplies these maps by an element of $G$ and/or $\pm 1$ and this
does \emph{not} affect $\fS_{K/k}$ at all.\vspace{1ex}\\
\example\label{ex:2A}\ {\bf The Case $k=\bbQ$.\ }
Suppose $k=\bbQ$ and $\fm$
is the cycle $f\bbZ.\infty$ for some integer $f>1$.
Thus $k(\fm)$ is the field $\bbQ(\mu_f)$. We may identify
$\Cl_{\fm}(\bbQ)$ with $(\bbZ/f\bbZ)^\times$ by $[(a)]_\fm\leftrightarrow \bar{a}$ for
integers $a>0$ with $(a,f)=1$. The Artin map then sends $\bar{a}$ to $\sigma_a\in G_{\fm}$
where $\sigma_a(\zeta)=\zeta^a$ for every $\zeta\in\mu_f$.
Working through the definitions in this case and setting $\xi_f=\be(1/f)$,
equation\refeq{2A.9} becomes
\[
\Phi_{\fm}(s)=\sum_{a=1\atop (a,f)=1}^{f}Z(s;a,f)\sigma_a\inv
\mbox{\ \ where\ \
$Z(s;a,f):={\displaystyle \sum_{n\geq 1}\frac{\xi_f^{an}}{n^s}}$ for ${\rm Re}(s)>1$}
\]
Standard methods of analytic continuation
(\eg\ a very simple case of~\cite[Prop.~1]{Shi}) give $Z(0;a,f)=\xi_f^a/(1-\xi_f^a)$.
For any complex abelian extension $K$ of $\bbQ$ with $\fm(K)=\fm$ it follows that
\beql{eq:2K}
\Phi_{K/\bbQ}(0)=\frac{1}{f}\sum_{a=1\atop (a,f)=1}^{f}
\left(\frac{\xi_f^a}{1-\xi_f^a}\right)\sigma_a\inv|_K
\eeq
To simplify, we now suppose that the prime number $p$ is odd and
we take $K=\bbQ(\mu_{p^{n+1}})$ for some
$n\geq 0$ in the above. Thus $\fm=\fm(K)=p^{n+1}\bbZ.\infty$
and $\bbQ(\fm)$ equals $K$ which has a unique, totally ramified prime $\fP=\fP_1$
above $p$. If we also take $\tau_1$ to be the identity then $j=j\tau_1$
defines an isomorphism from $K_\fP$ to $\overline{j(K)}=\bbQ_p(\mu_{p^{n+1}})$
and by total ramification we have an isomorphism
$G\rightarrow \Gal(\bbQ_p(\mu_{p^{n+1}})/\bbQ_p)$ sending $g$ to $\hat{g}$ such that
$jg(x)=\hat{g}j(x)$ for all $x\in K_\fP$. Thus if $u$ is any element of
$\bigwedge^1_{\bbZ_p G}U^1_p(K)=U^1_p(K)=U^1(K_\fP)$ then, in this notation,
\[
R_{K/\bbQ,p}^{(j)}(u)=\lambda^{(j)}_{1,p}(u)=
\sum_{b=1\atop (b,p)=1}^{p^{n+1}}\log_p(j\sigma_b(u))\sigma_b\inv=
\sum_{b=1\atop (b,p)=1}^{p^{n+1}}\hat{\sigma}_b(\log_p(j(u)))\sigma_b\inv
\]
On the other hand, if we let
$\zeta_n:=j(\xi_{p^{n+1}})$ then equation\refeq{2K} gives
\[
j(\sqrt{d_{\bbQ}}\Phi_{K/\bbQ}(0)^\ast)=
\frac{1}{p^{n+1}}\sum_{a=1\atop (a,p)=1}^{p^{n+1}}
\left(\frac{\zeta_n^a}{1-\zeta_n^a}\right)\sigma_a=
\frac{1}{p^{n+1}}\sum_{a=1\atop (a,p)=1}^{p^{n+1}}
\hat{\sigma}_a\left(\frac{\zeta_n}{1-\zeta_n}\right)\sigma_a
\]
If we write $\alpha_c$ for the coefficient of $\sigma_c\inv$ in $\fs_{K/\bbQ}(u)$ then
the last two equations imply that
\[
\alpha_c=\frac{1}{p^{n+1}}\Tr_{\bbQ_p(\mu_{p^{n+1}})/\bbQ_p}\left(\frac{\zeta_n}{1-\zeta_n}
\log_p(j\sigma_c(u))\right)
\]
where, of course, $j\sigma_c(u)$ lies in $U^1(\bbQ_p(\mu_{p^{n+1}}))$.
The right-hand side is familiar from the explicit reciprocity law proved by
Artin and Hasse in~\cite{AH} (see also~\cite[Thm.~10, Ch.~12]{AT} for the case $n=0$).
More precisely, their law states firstly (and implicitly) that
\beql{eq:firstfact}
\alpha_c\in\bbZ_p\ \ \ \mbox{for all $c$ and all $u\in U^1(K_\fP)$}
\eeq
which is not \emph{a priori} obvious, and secondly that
\beql{eq:secondfact}
\zeta_n^{\alpha_c}=((1-\zeta_n)\inv,j\sigma_c(u))_\pnpo\ \ \
\mbox{for all $c$ and all $u\in U^1(K_\fP)$}
\eeq
where $ (\cdot,\cdot)_\pnpo:\bbQ_p(\mu_\pnpo)^\times\times
\bbQ_p(\mu_\pnpo)^\times\rightarrow\mu_\pnpo$ is the Hilbert symbol with
values in $\mu_\pnpo$. (Thus $(\alpha,\beta)_\pnpo$ equals $(\frac{\alpha}{\beta})$)
in the notation of~\cite{AH}.)
Of course,\refeq{firstfact} is simply the statement that
$\fS_{K/\bbQ}$ is contained in $\bbZ_p G$.
This motivates the more general question of the $p$-integrality
of $\fS_{K/k}$ which we shall address in the next section
and in Subsection\ref{subsec:thm3arev}.
Equation\refeq{secondfact} amounts to a congruence for $\fs_{K/\bbQ}(u)$
modulo $\pnpo$ of which we shall give a conjectural generalisation
in Subsection\ref{subsec:rubinstark}.
\section{Integrality Properties of $\SKk$}\label{sec:integrality}
We shall investigate an integrality condition on $\SKk$ as
a $\bbZ_p G$-submodule of $\bbQ_p G$.
Largely for the sake of simplicity, \emph{we shall assume henceforth}
\begin{hyp}
$p$ is odd.
\end{hyp}
and we introduce the notation $\delta(K/k)=\delta_p(K/k):=
|\{i\,:\, 1\leq i\leq d,\fp_i\nin S_{\rm ram}(K/k)\}|$. A preliminary
result is
\begin{prop}\label{prop:3A}
If $e_{K/k}^-=0$ then $\SKk\subset p^{\delta(K/k)}\bbZ_p G$.
\end{prop}
If $k\neq \bbQ$ then Prop.~\ref{prop:2A} gives $\PhiKko=0$ and the result is trivial.
On the other hand,
$\Phi_{K/\bbQ}(0)=-\half\prod_{q|f(K)}(1-q\inv)e_{\chi_0}$.
The result follows from this formula, properties of $\log_p$ and the fact that
$N_{K/\bbQ}(U^1_p(K))=N_{\bbQ(\mu_{f(K)})/\bbQ}(U^1_p(\bbQ(\mu_{f(K)})))
\subset 1+f(K)\bbZ_p$. We omit the details
since Cor~\ref{cor:3A}
is much more general.\ePf
\noindent\
Equation~\refeq{introA} implies that $\Theta_{K/k}(0)$ has coefficients in
$w_K\inv\bbZ$. Equation\refeq{2J.5} (and the obvious fact that $w_{K[\fn]}|w_K$
for all $\fn$) therefore implies
\beql{eq:3A}
\PhiKko\in(w_Kd_kN\ff(K))\inv\bbZ[\mu_{f(K)}]G
\eeq
Now suppose that $p$ is unramified in $K/\bbQ$. Then
$p\ndiv w_Kd_k N\ff(K)$ (recall: $p\neq 2$) so equation\refeq{3A} implies
$j(\PhiKkoa)\in\bar{\bbZ}_p G$ where $\bar{\bbZ}_p$ denotes the
integral closure of $\bbZ_p$ in $\barbbQ_p$. Furthermore, in this case,
$\overline{j\tau_i(K)}$ is an unramified
extension of $\bbQ_p$ for each $i$. This implies
$\lambda_{i,p}(u)\in p\bar{\bbZ}_p G$ for all $i$ and $u\in U^1_p(K)$ and
hence $R^{(j)}_{K/k,p}(\theta)\in p^d\bar{\bbZ}_p G$
for all $\theta\in\bigwedge_{\bbZ G}^dU^1_p(K)$. From these considerations, we deduce easily
\begin{prop}\label{prop:unram}
If $p$ is unramified in $K/\bbQ$ then
$\SKk\subset p^{\delta(K/k)}\bbZ_p G$.\ePf
\end{prop}
Returning to the case in which $p$ may ramify in $K/\bbQ$,
the above argument makes it clear that we can expect non-$p$-integral
coefficients both in $\PhiKkoa$ (coming from $w_K\inv$, $d_k\inv$ and/or $N\ff(K)\inv$)
and in $R^{(j)}_{K/k,p}(\theta)$ (coming from $\log_p(x)$, $x\in U^1(\overline{j\tau_i(K)})$).
The proof of the following result is therefore much more delicate than the above
and occupies the rest of this section.
\begin{thm}\label{thm:3A}
Suppose that $p$ splits completely in $k$ and also that $K$ satisfies
at least one of the following two conditions:
\beql{eq:3B}
\mbox{there exists $\fq\in S_{\rm ram}(K/k)$ not dividing $p$, or}
\eeq
\beql{eq:3C}
p\ndiv w_K
\eeq
Then $\SKk\subset p^{\delta(K/k)}\bbZ_p G$.
\end{thm}
\begin{cor}\label{cor:3A}
If $k=\bbQ$ then
$\fS_{K/\bbQ}\subset p^{\delta(K/\bbQ)}\bbZ_p G$.\ePf
\end{cor}
\noindent \textsc{Proof of Corollary~\ref{cor:3A}\ } If Condition\refeq{3C}
fails then $\bbQ(\mu_p)\subset K$. In particular, $K$ is complex
so if Condition\refeq{3B} also
fails then $\fm(K)=p^{n+1}.\infty$ for some $n\geq 0$, hence
$K\subset\bbQ(\fm(K))=\bbQ(\mu_{p^{n+1}})$. The last two
inclusions -- and the minimality of $n$ in the second -- imply
$K=\bbQ(\mu_{p^{n+1}})$ and this case was dealt with explicitly
in Example~\ref{ex:2A}.\ePf
\noindent \textsc{Proof of Theorem~\ref{thm:3A}\ }
Since $\fp$ splits in $k$, the primes $\fp_1,\ldots,\fp_d$ are distinct
and we can uniquely write
\[
\fm(K)=\ff(K)\fz(K)=\ff'(K)\fp_1^{n_1(K)+1}\ldots\fp_1^{n_d(K)+1}\fz(K)
\]
where $\ff'(K)$ is prime to $p$ and $n_i(K)\geq -1$ for $i=1,\ldots,d$.
Condition\refeq{3B} holds if and only if  $\ff'(K)\neq\cO$.
If it fails, we proceed as follows.
Lemme~IV.1.1 of~\cite{Tate} implies that
\[
w_K={\rm gcd}\{N\fq-1\,:\,\fq\in \cQ \}
\]
where $\cQ$ denotes the set of of all prime ideals of $k$ such that
$\fq\ndiv pw_K$, $\fq\nin S_{\rm ram}(K/k)$ and $\sigma_{\fq,K/k}=1$.
By assumption, Condition\refeq{3C} must hold, so
we can find $\fq_0\in\cQ$ such that $p\ndiv N\fq_0-1$.
Now set
\[
\ff'=
\left\{
\begin{array}{ll}
\ff'(K)&
    \mbox{if Condition\refeq{3B} holds and}\\
\fq_0&
    \mbox{if not.}
\end{array}
\right.
\]
Thus $\ff'\neq\cO$ and $(\ff',p)=1$ in both cases and we define an ideal $\ff$ and
a cycle $\fm$ by
\beql{eq:3C.5}
\fm:=\ff\uinf=\ff'\fp_1^{n_1+1}\ldots\fp_1^{n_d+1}\uinf
\eeq
where $n_i={\rm max}\{n_i(K),0\}$ and $\uinf$ denotes the `product' of all the
infinite (real) places of $k$. Since $\fm(K)|\fm$, we have $k(\fm)\supset K$.
(In fact, it is not hard to see that
$\fm=\fm(K(\mu_p))$ or $\fq_0\fm(K(\mu_p))$ according as
Condition\refeq{3B} holds or fails. In the second
case, $\fm$ is not necessarily a conductor but in either case
we actually have $k(\fm)\supset K(\mu_p)$).
Theorem~\ref{thm:3A} is now a consequence of the two following statements.
\begin{claim}\label{claim:3A}
$k(\fm)/K$ is at most tamely ramified at all primes dividing $p$.
\end{claim}
\begin{claim}\label{claim:3B}
$j(N\ff\inv\Phi_\fm(0)^\ast)R^{(j)}_{k(\fm)/k,p}(\theta)$ lies in
$\bar{\bbZ}_p G_\fm$
for every $\theta\in\bigwedge_{\bbZ G_\fm}^dU^1_p(k(\fm))$.
\end{claim}
Let us see how Claims~\ref{claim:3A}
and~\ref{claim:3B} imply the Theorem. Using equation\refeq{2B} we find
\[
\pi_{k(\fm),K}(N\ff\inv \Phi_\fm(0)^\ast)=
\pi_{k(\fm(K)),K}\circ \pi_{\fm,\fm(K)}(N\ff\inv \Phi_\fm(0)^\ast)=
p^{-\delta_p(K/k)}A\sqrt{d_k}\PhiKkoa
\]
where $A$ is either $\sqrt{d_k}\prod_{\fp|p,\,\fp\ndiv\ff(K)}(p-\sigma_{\fp,K/k})$
or this element multiplied by $1-N\fq_0\inv$ according as
Condition\refeq{3B} holds or fails. In any case it is easy to see that
$j(A)$ lies in $\bbZ_p G^\times$. (The fact that
$p$ splits in $k$ implies $j(\sqrt{d_k})\in\bbZ_p^\times$.)
Now Claim~\ref{claim:3A} implies that the norm $N_K:=N_{k(\fm)/K}$ defines a
surjective map from $U^1_p(k(\fm))$ to $U^1_p(K)$ and hence a unique surjection
from $\bigwedge_{\bbZ G_\fm}^d U^1_p(k(\fm))$ to $\bigwedge_{\bbZ G}^d U^1_p(K)$
-- also denoted $N_K$ --
sending $u_1\wedge\ldots\wedge u_d$ to $N_Ku_1\wedge\ldots\wedge N_Ku_d$.
It is easy to check that
$R^{(j)}_{K/k,p}\circ N_K=\pi_{k(\fm),K}\circ R^{(j)}_{k(\fm)/k,p}$.
Thus given any $\phi\in\bigwedge_{\bbZ G}^d U^1_p(K)$ we may choose
$\theta\in\bigwedge_{\bbZ G_\fm}^d U^1_p(k(\fm))$ such that $\phi=N_K(\theta)$ and we find
\begin{eqnarray*}
p^{-\delta_p(K/k)}\sKk(\phi)&=&p^{-\delta_p(K/k)}j(\sqrt{d_k}\PhiKkoa)R^{(j)}_{K/k,p}(\phi)\\
&=&j(A)\inv\pi_{k(\fm),K}(j(N\ff\inv \Phi_\fm(0)^\ast)
R^{(j)}_{k(\fm)/k,p}(\theta))
\end{eqnarray*}
which lies in $\bar{\bbZ}_p G$ by Claim~\ref{claim:3B} and also in $\bbQ_p G$
hence in $\bbZ_p G$. This gives the Theorem.

Some preliminary notation and lemmas will establish
Claim~\ref{claim:3A} and pave the way for the proof of the
harder Claim~\ref{claim:3B}. Suppose that $L/F$ is a finite,
abelian field extension. If $L$ and $F$ are local fields and
$l\geq -1$ is an integer, we write $G(L/F)^l$ for the $l$th ramification
group in the `upper numbering' and if $L$ and $F$ are number
fields and $\fq$  a prime of $\cO_F$, then $G(L/F)^l_\fq$ will
denote the $l$th ramification group at any prime $\fQ$ of $\cO_L$
above $\fq$, naturally identified with $G(L_\fQ/F_\fq)^l$. Taking
$F$ to be $k$, local and global class field theory give the
formula \beql{eq:3D} \ord_\fq(\ff(L))=
\min\{l\in\bbN:G(L/k)^l_\fq=\{1\}\}\ \ \ \mbox{for all primes
$\fq$ of $\cO$} \eeq
For any integer $r\geq -1$ we shall use the
abbreviation $\mu(r)$ for $\mu_{p^{r+1}}$, either as a
subgroup of $\barbbQ^\times$ or of $\barbbQ_p^\times$. Similarly,
$\zeta_r$ will denote either $\be(1/p^{r+1})$ or its embedding in $\barbbQ_p$ under $j$.
It is well known that, for any $l\geq 0$
\beql{eq:3E}
\mbox{$G(\bbQ_p(\zeta_r)/\bbQ_p)^l$ equals
$\Gal(\bbQ_p(\zeta_r)/\bbQ_p(\zeta_{l-1}))$
if $l\leq r$, $\{1\}$ if not}
\eeq
(see \cite[Ch.~IV]{SerreLF}). It follows easily from the above
that $\fm(k(\zeta_r))=p^{r+1}\cO\uinf$.
Given any $i=1,\ldots,d$ and a cycle $\fn=\fg\fz$
we shall write $\fn^{(i)}$ for its prime-to-$\fp_i$ part
so that $\fn=\fp_i^{r_i+1}\fn^{(i)}$ where $r_i+1=\ord_{\fp_i}(\fg)$.
We also write $D_{\fn,i}$ for $D_{\fp_i}(k(\fn)/k)=G(k(\fn)/k)^{-1}_{\fp_i}$
and $T_{\fn,i}$ for $T_{\fp_i}(k(\fn)/k)=G(k(\fn)/k)^{0}_{\fp_i}$ so that
$T_{\fn,i}\subset D_{\fn,i}\subset G_\fn$ and
$T_{\fn,i}=\Gal(k(\fn)/k(\fn^{(i)}))=\ker(\pi_{\fn,\fn^{(i)}})$.
By mapping the sequence\refeq{2I} for $\fm=\fn$ onto the one for $\fn^{(i)}$, we see that
$T_{\fn,i}$ is a quotient of $(\cO/\fp_i^{r_i+1})^\times$ and since
$\fp_i$ splits in $k$ we obtain
\begin{lemma}\label{lemma:3A}
The inertia group $T_{\fn,i}$ is a quotient of $(\bbZ/p^{r_i+1}\bbZ)^\times$
where $r_i+1=\ord_{\fp_i}(\fn)$.\ePf
\end{lemma}
We also set $f_{\fn,i}:=|D_{\fn,i}:T_{\fn,i}|=|D_{\fn^{(i)},i}|$
(the residual degree of $k(\fn)/k$ above $\fp_i$).
Now let $n={\rm max}\{n_1,\ldots,n_d\}={\rm max}\{0,n_1(K),\ldots,n_d(K)\}$
and let $\tfm$ be the cycle $\ff'p^{n+1}\cO\uinf$, so that $\fm(K)|\fm|\tfm$
and $K\subset k(\fm(K))\subset k(\fm)\subset k(\tfm)$.
\begin{lemma}\label{lemma:3A.5}
For each $i=1,\ldots,d$ we have $k(\tfm)=k(\tfm^{(i)})k(\zeta_n)$ and
$k(\tfm^{(i)})\cap k(\zeta_n)=k$ so that $\Gal(k(\tfm)/k)=T_{\tfm,i}\times
\Gal(k(\tfm)/k(\zeta_n))$
and the restriction maps $\Gal(k(\tfm)/k(\zeta_n))\rightarrow G_\tfmi$
and $T_{\tfm,i}\rightarrow\Gal(k(\zeta_n)/k)$ are isomorphisms.
\end{lemma}
\bPf\ Since $\fm(k(\zeta_n))=p^{n+1}\cO\uinf$ divides $\tfm$ it follows
that $k(\zeta_n)\subset k(\tfm)$ and
$k(\zeta_n)/k$ is totally ramified at each $\fp_i$ since $p$ splits in $k$.
Hence the reduction map $T_{\tfm,i}\rightarrow\Gal(k(\zeta_n)/k)$ is onto.
But the previous lemma shows that $T_{\tfm,i}$ is a quotient of
$(\bbZ/p^{n+1}\bbZ)^\times\cong\Gal(k(\zeta_n)/k)$ so this map must be an
isomorphism. We have already observed that $T_{\tfm,i}=\Gal(k(\tfm)/k(\tfmi))$
so the result follows.\ePf
\noindent The lemma shows that the action of $\Ttmi$ on $\mu(n)$ defines an
isomorphism onto $(\bbZ/p^{n+1}\bbZ)^\times$. Its inverse gives rise to a surjection
$\tsigma_i:\bbZ_p^\times\rightarrow T_{\tfm,i}$ fitting into the exact sequence
\[
1\rightarrow (1+p^{n+1}\bbZ_p)\longrightarrow \bbZ_p^\times
\stackrel{\tsigma_i}{\longrightarrow}
\Ttmi\rightarrow 1
\]
The isomorphism $T_{\fm,i}\rightarrow\Gal(k(\zeta_n)/k)$ respects higher
ramification groups in the upper numbering,
so\refeq{3E} implies
\begin{lemma}\label{lemma:3B}
$\tsigma_i(1+p^l\bbZ_p)=G(k(\tfm)/k)^l_{\fp_i}$ for all $l\geq 1$ and $i=1,\ldots,d$.
In particular, if $1\leq l\leq n$ then the unique maximal subgroup of
$G(k(\tfm)/k)^l_{\fp_i}$ equals $G(k(\tfm)/k)^{l+1}_{\fp_i}$ and is of index $p$.
\ePf
\end{lemma}
\begin{lemma}\label{lemma:3C}
Suppose $n_i(K)\geq 0$ for some $i$ so that $n_i=n_i(K)$ (by definition). Then
\beql{eq:3F}
G(k(\tfm)/K)^1_{\fp_i}=G(k(\tfm)/k(\fm(K)))^1_{\fp_i}=G(k(\tfm)/k(\fm))^1_{\fp_i}=
G(k(\tfm)/k)^{n_i+1}_{\fp_i}
\eeq
\end{lemma}
\bPf\ The inclusions $K\subset k(\fm(K))\subset k(\tfm)$ give
$G(k(\tfm)/K)^1_{\fp_i}\supset G(k(\tfm)/k(\fm(K)))^1_{\fp_i}\supset
G(k(\tfm)/k(\fm))^1_{\fp_i}$ (take Sylow $p$-subgroups of inertia).
Lemma~\ref{lemma:3A} shows that $G(k(\tfm)/k(\fm))^1_{\fp_i}$ is a
$p$-subgroup of index at most $p^{n_i}$ in $G(k(\tfm)/k)^{1}_{\fp_i}$,
so Lemma~\ref{lemma:3B} gives the inclusion
$G(k(\tfm)/k(\fm))^1_{\fp_i}\supset G(k(\tfm)/k)^{n_i+1}_{\fp_i}$ (which also follows from
$G(k(\fm)/k)^{n_i+1}_{\fp_i}=\{1\}$.) If the resulting inclusion
$G(k(\tfm)/K)^1_{\fp_i}\supset G(k(\tfm)/k)^{n_i+1}_{\fp_i}$ were strict then
we would have ($n_i\geq 1$ and)
$G(k(\tfm)/K)^1_{\fp_i}\supset G(k(\tfm)/k)^{n_i}_{\fp_i}$, by Lemma~\ref{lemma:3B} again.
This would imply $G(K/k)^{n_i}_{\fp_i}=1$ so that
$n_i(K)+1=\ord_{\fp_i}(\fm(K))$ could be at most $n_i$,
by\refeq{3D}, contradicting our assumption. We conclude that
$G(k(\tfm)/K)^1_{\fp_i}=G(k(\tfm)/k)^{n_i+1}_{\fp_i}$ and the result follows.\ePf
\noindent
{\sc Proof of Claim~\ref{claim:3A}}\
The equality of the first and the third groups in\refeq{3F} show that
if $n_i(K)\geq 0$ then the extension $k(\fm)/K$ is at most tamely ramified above $\fp_i$.
On the other hand, if $n_i(K)=-1$ then $n_i=0$ so the same holds of the whole extension
$k(\fm)/k$ (\eg\ by Lemma~\ref{lemma:3A} with $\fn=\fm$).\ePf
\noindent From now on we can forget the fields $K$ and $k(\fm(K))$ and concentrate
on $k(\fm)$ and $k(\tfm)$. Since $\fm(k(\zeta_p))=p\cO.\uinf$
which divides $\fm$, we have $k(\zeta_p)\subset k(\fm)$. In particular, $p-1$ divides
$e_{\fp_i}(k(\fm)/k)$ for all $i$ so
Lemma~\ref{lemma:3A} implies that the tame ramification
indices above $p$ in $k(\tfm)/k(\fm)$ are all $1$. Thus
\[
\ker(\pi_{\tfm,\fm}:\Ttmi\rightarrow\Tmi)=
G(k(\tfm)/k(\fm))^0_{\fp_i}=G(k(\tfm)/k(\fm))^1_{\fp_i}=
G(k(\tfm)/k)^{n_i+1}_{\fp_i}
\]
Hence by Lemma~\ref{lemma:3B} we may deduce the exactness of the sequence
\beql{eq:3G}
1\rightarrow (1+p^{n_i+1}\bbZ_p)\longrightarrow \bbZ_p^\times
\stackrel{\pi_{\tfm,\fm}\circ\tsigma_i}{\longrightarrow} \Tmi\rightarrow 1
\eeq
Now
Lemma~\ref{lemma:3A.5} implies that for each $i=1,\ldots,d$ there exists a unique element
$\tphi_i$ of $D_{\tfm,i}$ which restricts to
the Frobenius $\sigma_{\fp_i,\tfm^{(i)}}\in D_{\tfm^{(i)},i}$ on $k(\fm^{(i)})$
and fixes $\mu(n)$.
\begin{lemma}\label{lemma:3C.5}
Fix $i\in\{1,\ldots,d\}$, let $Z$ be any set of representatives
of $\bbZ_p^\times$ modulo $1+p^{n_i+1}\bbZ_p$ and $L$ any set of
representatives of  $\bbZ$ modulo $f_{\fm,i}\bbZ$.
Then the map $Z\times L\rightarrow D_{\fm,i}$ sending
$(z,l)$ to $\pi_{\tfm,\fm}(\tsigma_i(z)\tphi_i^l)$ is bijective.
\end{lemma}
\bPf\ The exact sequence\refeq{3G} shows that
$|Z\times L|=|\Tmi|f_{\fm,i}=|D_{\fm,i}|$ so it suffices to prove injectivity.
But if $\pi_{\tfm,\fm}(\tsigma_i(z')\tphi_i^{l'})=\pi_{\tfm,\fm}(\tsigma_i(z)\tphi_i^l)$ then
applying $\pi_{\fm,\fm^{(i)}}$, we find $\sigma_{\fp_i,\tfm^{(i)}}^{l'}=\sigma_{\fp_i,\tfm^{(i)}}^l$
so $l'\equiv l\pmod{f_{\fm,i}}$ so $l=l'$ and
$\pi_{\tfm,\fm}(\tsigma_i(z'))=\pi_{\tfm,\fm}(\tsigma_i(z))$
which implies $z'=z$ by\refeq{3G} again.\ePf
\noindent We next introduce some notation and preparatory lemmas
in the \emph{local} situation. We shall write $\fP_{\tfm,i}$ and
$\fP_{\fm,i}$ for those primes above $\fp_i$ in $k(\tfm)$ and $k(\fm)$ respectively
which are induced by $j\tau_i$. We set
$M_{\tfm,i}:=\overline{j\tau_i(k(\tfm))}$
which contains $M_{\fm,i}:=\overline{j\tau_i(k(\fm))}$,
isomorphic via $j\tau_i$ to the completions $k(\tfm)_{\fP_{\tfm,i}}$ and
$k(\fm)_{\fP_{\fm,i}}$ respectively.
Lemma~\ref{lemma:3A.5} and the splitting hypothesis show that
$M_{\tfm,i}=H_{f_{\tfm,i}}(\zeta_n)$
and where $H_{f}$ denotes the unique unramified extension
of $\bbQ_p$ of degree $f$ inside $\barbbQ_p$.
The group $D_{\tfm,i}$ is isomorphic to $\Gal(M_{\tfm,i}/\bbQ_p)=
\Gal(M_{\tfm,i}/H_{f_{\tfm,i}})\times \Gal(M_{\tfm,i}/\bbQ_p(\zeta_n))$
by the map $g\mapsto \check{g}$ where $\check{g}\circ(j\tau_i)=(j\tau_i)\circ g$
on $k(\tfm)$.
Now choose an unramified extension $H$ of $\bbQ_p$ which contains $H_{f_{\tfm,i}}$
for $i=1,\ldots,d$. Let $\mu(\infty):=\bigcup_{m\geq 0}\mu(n)$ be
the group of all $p$-power roots of unity and write
$H(\zeta_\infty)$ for $H(\mu(\infty))$ \etc\ Thus
$\Gal(H(\zeta_\infty)/\bbQ_p)=
\Gal(H(\zeta_\infty)/H)\times \Gal(H(\zeta_\infty)/\bbQ_p(\zeta_\infty))$.
The action on $\mu(\infty)$ defines an isomorphism
$\hat{\sigma}:\bbZ_p^\times\rightarrow\Gal(H(\zeta_\infty)/H)$ and
shall write $\hat{\phi}$ for the unique lift to $\Gal(H(\zeta_\infty)/\bbQ_p(\zeta_\infty))$
of the local Frobenius element $\phi\in\Gal(H/\bbQ_p)$.
These notations clearly link up with the global ones as follows
\beql{eq:3H}
\mbox{
$\hat{\sigma}(z)|_{M_{\tfm,i}}=\tsigma_i(z)\check{}$
and
$\hat{\phi}|_{M_{\tfm,i}}=\check{\tphi_i}$
for all $z\in\bbZ_p^\times$ and $i=1,\ldots,d$}
\eeq
\begin{lemma}\label{lemma:3D} For each $i$ we have inclusions
$M_{\fm,i}\subset H_{f_{\tfm,i}}(\zeta_{n_i})\subset H(\zeta_{n_i})$. Moreover
the extension $H(\zeta_{n_i})/M_{\fm,i}$ is unramified.
\end{lemma}
\bPf Since $e(H(\zeta_{n_i})/\bbQ_p)=(p-1)p^{n_i}=e(M_{\fm,i}/\bbQ_p)$ by
the exact sequence\refeq{3G}, it suffices to prove the first of the two inclusions.
Now $M_{\fm,i}\subset M_{\tfm,i}=H_{f_{\tfm,i}}(\zeta_n)$ so let
$\sigma$ be any element of $\Gal(H_{f_{\tfm,i}}(\zeta_n)/H_{f_{\tfm,i}}(\zeta_{n_i}))$.
Clearly, $\sigma$ may be lifted to $\hat{\sigma}(z)\in\Gal(H(\zeta_\infty)/H)$ for
some $z\in 1+p^{n_i+1}\bbZ_p$. The sequence\refeq{3G} shows that
$\tsigma_i(z)$ fixes $k(\fm)$ and it follows from\refeq{3H} that
$\hat{\sigma}(z)$ fixes $M_{\fm,i}$ \ie\
$\sigma\in \Gal(H_{f_{\tfm,i}}(\zeta_n)/M_{\fm,i})$. The result follows\ePf
\rem\ There is in general
no reason to suppose that $M_{\fm,i}$ is contained in
$H_{f_{\fm,i}}(\zeta_{n_i})$, even though
these two fields have the same residual degree over $\bbQ_p$. Since their absolute
ramification indices are also equal, such an inclusion would actually imply equality
and occurs if and only if $\zeta_{n_i}\in M_{\fm,i}$.\vspace{1ex}\\
We define a unique action (denoted `$\,\cdot\,$')
of $\Gal(H(\zeta_\infty)/\bbQ_p)$ on the power series
ring $H[[X]]$ by setting
\[
\mbox{
$\hat{\sigma}(z)\cdot F(X)=F((1+X)^z-1)$
and
$\hat{\phi}\cdot F(X)=\hat{\phi}(F(X))$
for all $z\in\bbZ_p^\times$ and $F\in H[[X]]$}
\]
(where $\hat{\phi}(\sum_{m\geq 0}a_iX^i)$ means $\sum_{m\geq 0}\hat{\phi}(a_i)X^i$).
For any $\sigma\in \Gal(H(\zeta_\infty)/\bbQ_p)$ the
map $F\mapsto \sigma\cdot F$ is a ring automorphism
preserving $\cO_H[[X]]$.
Moreover, if $F$ lies in $\cO_H[[X]]$ then $F(\zeta-1)$ converges
for all $\zeta\in \mu(\infty)$ and $\sigma(F(\zeta-1))=(\sigma\cdot F)(\zeta-1)$
(which also equals $F(\sigma(\zeta-1))$ if $\sigma\in \Gal(H(\zeta_\infty)/H)$).
The following result, which is a key ingredient in the proof of Claim~\ref{claim:3B}, has
its origins in Coleman's proof of~\cite[Thm.~26]{Coleman1}.
\begin{lemma}\label{lemma:3E} Suppose $i\in\{1,\ldots,d\}$ and $\hat{u}\in U^1(\Mmi)$ are given.
Then there exists a power series $g\in H[[X]]$ with the following properties
\begin{enumerate}
\item\label{A} $g\in\cO_H[[X]]$,
\item\label{B} $g(0)=0$,
\item\label{D} $g(\zeta-1)\in\Mmi\ \forall\,\zeta\in\mu(n_i)$ and
\item\label{C} ${\displaystyle\log_p(\hat{u})=\sum_{t=0}^{n_i}
\frac{\hat{\phi}^t(g(\zeta_{n_i-t}-1))}{p^t}}$.
\end{enumerate}
\end{lemma}
\bPf\ Since $H(\zeta_{n_i})/M_{\fm,i}$ is unramified
by Lemma~\ref{lemma:3D}, we may choose
$\tilde{u}\in U^1(H(\zeta_{n_i}))$ such that
$\hat{u}=N_{H(\zeta_{n_i})/M_{\fm,i}}(\tilde{u})$. Since $\zeta_{n_i}-1$ is a uniformizer
for $H(\zeta_{n_i})$, we may clearly choose $\tilde{h}\in 1+X\cO_H[[X]]$
such that $\tilde{u}=\tilde{h}(\zeta_{n_i}-1)$. Choose also a
set $\Sigma\subset \Gal(H(\zeta_\infty)/\bbQ_p)$ of representatives modulo
$\Gal(H(\zeta_\infty)/\Mmi)$ and set
$h(X)=\prod_{\sigma\in\Sigma}(\sigma\cdot\tilde{h}(X))$. Thus $h\in 1+X\cO_H[[X]]$ and
\beql{eq:3J}
h(\zeta-1)=N_{H(\zeta_{n_i})/M_{\fm,i}}\tilde{h}(\zeta-1)
\in U^1(\Mmi)\ \forall\,\zeta\in\mu(n_i)
\eeq
In particular, $h(\zeta_{n_i}-1)=\hat{u}$. Let
\[
g=\log(h(X))-\frac{\hat{\phi}(\log(h(1+X)^p-1))}{p}=\frac{1}{p}
\log(h(X)^p/\hat{\phi}(h((1+X)^p-1)))
\]
where $\log:1+XH[[X]]\rightarrow H[[X]]$ is the homomorphism
sending the power series $a(X)$ to $\sum_{m=1}^\infty(-1)^{m-1}(a(X)-1)^m/m$.
Since $h(X)^p$ and $\hat{\phi}(h((1+X)^p-1))$ lie in the multiplicative group
$1+X\cO_H[[X]]$, so does their quotient and reducing modulo $p$ it is easy to see that
the latter actually lies in $1+pX\cO_H[[X]]$. Applying $\log$, it follows by
standard estimates that $g\in\cO_H[[X]]$, proving property~\ref{A}. Property~\ref{B}
is immediate.
For any $\zeta\in\mu(n_i)$ we have
\beql{eq:3K}
g(\zeta-1)=\log_p(h(\zeta-1))-\frac{\hat{\phi}(\log_p(h(\zeta^p-1)))}{p}
\eeq
(Identities such
as $\log(h)(\zeta_{n_i}-1)=\log_p(h(\zeta_{n_i}-1))$ \etc~may be rigorously treated
using \eg\ the one-variable analogue of~\cite[Prop.~5.3]{p-adic lims}).
Property~\ref{D}
follows by combining\refeq{3K} and\refeq{3J} as does Property~\ref{C} after a short computation,
noting that $h(0)=1$.\ePf
\noindent\ {\sc Proof of Claim~\ref{claim:3B}}
First note that $U^1_p(k(\fm))$ is generated as a $\bbZ G_\fm$-module
by the subgroups $U^1(k(\fm)_{\fP_{\fm,l}})$ for $l=1,\ldots,d$ so we only need
to treat the case $\theta=u_1\wedge\ldots\wedge u_d$ where each $u_i$ lies in
$U^1(k(\fm)_{\fP_{\fm,s(i)}})$ for some map $s:\{1,\ldots,d\}\rightarrow\{1,\ldots,d\}$.
Moreover
$\lambda_i$ vanishes on $U^1(k(\fm)_{\fP})$ for $\fP|p$ unless $\fP=\fP_{\fm,i}$.
It follows from the splitting hypothesis
that
the $\fP_{\fm,i}$ are distinct, hence that
$R^{(j)}_{k(\fm)/k,p}(\theta)=0$ unless $s$ is surjective, hence a permutation.
Thus we may assume w.l.o.g. that $u_i\in U^1(k(\fm)_{\fP_{\fm,i}})$. For the rest of the proof, we
fix such a choice of $u_i$, $i=1,\ldots,d$. It is easy to see that
$R^{(j)}_{k(\fm)/k,p}(u_1\wedge\ldots u_d)=\lambda_{1,p}(u_1)\ldots\lambda_{d,p}(u_d)$
and, moreover, that for each $i$
\[
\lambda_{i,p}(u_i)=
\sum_{g\in G_\fm}\log_p(j\tau_i(gu_i))g\inv=
\sum_{\delta\in D_{\fm,i}}\log_p(j\tau_i(\delta u_i))\delta \inv=
\sum_{\delta\in D_{\fm,i}}\check{\tilde{\delta}}(\log_p(\hat{u}_i)))\delta\inv
\]
where we set $\hat{u}_i=j\tau_i(u_i)\in U^1(M_{\fm,i})$
and $\tilde{\delta}$ is any lift of $\delta$ to $D_{\tfm,i}$
Furthermore, $N\ff\inv=N\ff'p^{-N}$ where $N:=n_1+\ldots n_d+d$ and
$N\ff'\in\bbZ_p^\times$. So, if we define $a_h$ in
$\barbbQ_p$ (actually in $\bbQ_p$)
by
\beql{eq:3L}
p^{-N}j(\Phi_\fm(0)^\ast)
\sum_{\delta_1\in D_{\fm,1}}\check{\tilde{\delta}}_1(\log_p(\hat{u}_1)))\delta_1\inv
\ldots\sum_{\delta_d\in D_{\fm,d}}\check{\tilde{\delta}}_d(\log_p(\hat{u}_d))\delta_d\inv=
\sum_{g\in G_\fm}a_hh
\eeq
then it suffices to show that $a_h\in \bar{\bbZ}_p$ for all $h\in G_\fm$.
Now the definition\refeq{2A.9} gives
$
j(\Phi_\fm(0)^\ast)=\sum_{g\in G_\fm}j(Z(0;g\cdot\fw_\fm^0))g
$
where we let $G_\fm$ act on $\fW_\fm$ via the Artin isomorphism
with $\Cl_\fm(k)$ (see below). Furthermore, for each $i$ we can choose
a power series
$g_i\in H[[X]]$ satisfying the four properties of Lemma~\ref{lemma:3E}
with respect to $\hat{u}=\hat{u}_i$. Substituting for $\log_p(\hat{u}_i)$
into\refeq{3L} and expanding gives
\beql{eq:3M}
a_h=p^{-N}\sum_{\delta_i\in D_{\fm,i}\atop i=1,\ldots,d}
\left(
j(Z(0;\delta_1\ldots\delta_dh\cdot\fw_\fm^0))
\prod_{i=1}^d\check{\tilde{\delta}}_i
\left(\sum_{t_i=0}^{n_i}
\frac{\hat{\phi}^{t_i}(g_i(\zeta_{n_i-t_i}-1))}{p^{t_i}}\right)
\right)=
p^{-N}\sum_{\underline{0}\leq \underline{t}\leq \underline{n}}a_{h,\underline{t}}
\eeq
say, where, in an obvious notation, the last sum ranges over all
$\underline{t}=(t_1,\ldots,t_d)$ such that $0\leq t_i\leq n_i\ \forall\,i$ and
\[
a_{h,\underline{t}}:=p^{-t_1-\ldots-t_d}
\sum_{\delta_i\in D_{\fm,i}\atop i=1,\ldots,d}
\left(
j(Z(0;\delta_1\ldots\delta_dh\cdot\fw_\fm^0))
\prod_{i=1}^d\check{\tilde{\delta}}_i\hat{\phi}^{t_i}(g_i(\zeta_{n_i-t_i}-1))
\right)
\]
(which makes sense, by Lemma~\ref{lemma:3E}~\ref{D}). For
each $i=1,\ldots,d$ we fix once and for all a set
$Z_i$ of representatives of $\bbZ_p^\times$ modulo
$1+p^{n_i+1}\bbZ_p$ and a set
$L_i$ of representatives of $\bbZ$ modulo
$f_{\fm,i}\bbZ$.
If $t$ is any integer then, as $l_i$ runs through $L_i$, so
$l_i-t$ runs through another set of representatives of $\bbZ$ modulo
$f_{\fm,i}\bbZ$. Hence, by Lemma~\ref{lemma:3C.5}, if $t_i$ is fixed then each
element $\delta_i$ of $D_{\fm,i}$ maybe written $\pi_{\tfm,\fm}(\tilde{\sigma}_i(z_i)
\tilde{\phi}_i^{l_i-t_i})$ for some unique $z_i\in Z_i$ and $l_i\in L_i$.
Thus, taking $\tilde{\delta}_i$ to be $\tilde{\sigma}_i(z_i)
\tilde{\phi}_i^{l_i-t_i}$ and applying\refeq{3H} we get
\begin{eqnarray}
\lefteqn{
a_{h,\underline{t}}:=p^{-t_1-\ldots-t_d}
\sum_{z_i\in Z_i\atop i=1,\ldots,d}
\sum_{l_i\in L_i\atop i=1,\ldots,d}
\left(
j(Z(0;\left(\prod_{i=1}^d\pi_{\tfm,\fm}(\tsigma_i(z_i)\tphi_i^{l_i-t_i})\right)
h\cdot\fw_\fm^0))\times\right.}
\hspace*{15em}&&\nonumber\\
 &&\left.\prod_{i=1}^d
\hat{\sigma}_i(z_i)\hat{\phi}^{l_i}
(g_i(\zeta_{n_i-t_i}-1))
\right)
\label{eq:3M.3}
\end{eqnarray}
Our aim is to rewrite this expression
for $a_{h,\underline{t}}$ in terms of certain $d$-variable
power series with coefficients in a finite extension of $H$.
This requires an interlude during which we recall some facts about $\fW_\fm$
and give some more auxiliary results:
An element of $\fW_\fm$ is an equivalence class
$\{\xi,I\}_\fm$ where the pair $(\xi,I)$ consists of a fractional
ideal $I$ of $k$ and a character of finite order $\xi:I\rightarrow\bbC^\times$
such that ${\rm ann}_\cO(\xi):=\{x\in\cO:\xi(xI)=0\}$ equals $\ff$.
(See~\cite{twizas} for more details here and in the following).
The equivalence relation is given (for $\fz=\uinf$) by
\beql{eq:3M.5}
\mbox{
$\{\xi,I\}_\fm=\{\xi',I'\}_\fm\Longleftrightarrow
\exists\, c\in k^\times_{\uinf}$ such that $I=cI'$ and
$\xi(ca')=\xi'(a')\ \forall\,a'\in I'$}
\eeq
(recall that $k^\times_{\uinf}$ is the group of totally positive
elements of $k$).
Suppose $\fw\in \fW_\fm$ is represented by  $(\xi,I)$
(\ie\ $\fw=\{\xi,I\}_\fm$) and $g=\sigma_{\fa,\fm}\in G_\fm$
for some integral ideal $\fa$ of $\cO$ prime to $\ff$. Then our definitions
and those of~\cite{twizas} give
\[
g\cdot\fw=[\fa]_\fm\cdot\fw:=\{\xi|_{\fa I},\fa I\}_\fm
\]
One checks that this depends only on the class $[\fa]_\fm\in\Cl_\fm(k)$
and gives a well-defined action
of $G_\fm$, via~$\Clmk$, on $\fW_\fm$. This action turns out to be free
and transitive; moreover $\fW_\fm$ contains a distinguished element $\fw_\fm^0$
(see~\cite{twizas}).

Now suppose that we are given a fractional ideal $J$ of $k$.
For simplicity we shall assume that $J$ is prime to $p$ which incurs no real loss
of generality in the sequel. Let $\rho:J\rightarrow\bbC^\times$
be a character of finite order
with ${\rm ann}_\cO(\rho)=\ff'$. This means that the image of $\rho$ is precisely
$\mu_{f'}$ where, as previously, $f'\in\bbZ$ denotes the
positive generator of $\ff'\cap \bbZ$. But $p\ndiv f'$ so
$\mu_{f'}$ is uniquely $p$-divisible and we may extend
$\rho$ uniquely to a character
\[
\rho:\bbZ[{\textstyle\frac{1}{p}}]J=
\bigcup_{\um\in\bbZ^d}\fp^\um J\longrightarrow\mu_{f'}\subset\bbC^\times
\]
where, for any $\um=(m_1,\ldots,m_d)\in\bbZ^d$,  the shorthand `$\fp^\um$' denotes
the fractional ideal $\fp_1^{m_1}\ldots\fp_d^{m_d}$. For any such $\um$ and
any $\us=(s_1,\ldots,s_d)\in\bbZ_p^d$, we define a character
$\xi_{\rho,\us,\um}\,:\,\fp^\um J\rightarrow \bbC^\times$
by
\beql{eq:3M.75}
\xi_{\rho,\us,\um}(a):=\rho(a)\zeta_n^{\sum_{i=1}^d s_ip^{-m_i}j\tau_i(a)}
\eeq
(The assumptions $a\in\fp^\um J$ and $(J,p)=1$ imply
$p^{-m_i}j\tau_i(a)\in\bbZ_p$ for all $i$).
Let $\cS(\un)$ denote the set
$\{\us\in\bbZ_p^d\,:\,\ord_p(s_i)=n-n_i, \forall\,i\}$.
By writing characters as the products of their $p$- and prime-to-$p$-part,
the reader should have no difficulty in proving the following
\begin{lemma}\label{lemma:3F} Suppose $J$ and $\um\in\bbZ^d$ are as above.
For any $\rho: J\rightarrow\bbC^\times$ with ${\rm ann}_\cO(\rho)=\ff'$
and $\us\in\cS(\un)$ we have ${\rm ann}_\cO(\xi_{\rho,\us,\um})=\ff$.
Conversely, if $\xi\,:\,\fp^\um J\rightarrow\bbC^\times$ is any
character with ${\rm ann}_\cO(\xi)=\ff$ then $\xi=\xi_{\rho,\us,\um}$
for some $\rho: J\rightarrow\bbC^\times$ with ${\rm ann}_\cO(\rho)=\ff'$
and $\us\in\cS(\un)$.
Moreover, $\rho$ is unique and $\us$ is unique modulo $p^n\bbZ_p^d$.\ePf
\end{lemma}
Given $\rho$, $J$, $\um$ and $\us$ as in the Lemma, we write $\fw_{\rho,\us,\um}$
for the class $\{\xi_{\rho,\us,\um},\fp^\um J\}_\fm\in\fW_\fm$.
\begin{lemma}\label{lemma:3G}  Suppose given
$\rho: J\rightarrow\bbC^\times$ with ${\rm ann}_\cO(\rho)=\ff'$,
$\um\in\bbZ^d$ and $\us\in\cS(\un)$. Then for any $\ul=(l_1,\ldots,l_d)\in\bbZ^d$
and $\uz=(z_1,\ldots,z_d)\in(\bbZ_p^\times)^d$, we have
\[
\prod_{i=1}^d\pi_{\tfm,\fm}(\tsigma_i(z_i)\tphi_i^{l_i})\cdot\fw_{\rho,\us,\um}=
\fw_{\rho,\uz\us,\ul+\um}
\]
where $\uz\us$ denotes $(z_1s_1,\ldots,z_ds_d)\in\cS(\un)$.
\end{lemma}
\bPf\ Since `$\,\cdot\,$' is an action, it
suffices to treat the case $\ul=(0,\ldots,0,1,0,\ldots,0)$
(`$1$' in the $i$th position) and $\uz=(1,\ldots,1,z,1,\ldots,1)$
($z\in\bbZ_p^\times$ in the $i$th position) for some $i$.
To simplify notation, we shall further assume $i=1$ (the other
$i$ being treated identically).
In this case, we choose
an integral ideal $\fa$ of $k$ prime to $\ff$
such that $\sigma_{\fa,\tfm}=\tsigma_1(z)\tphi_1\in D_{\tfm,1}\subset G_{\tfm}$
so that $\pi_{\tfm,\fm}(\tsigma_1(z)\tphi_1)=\sigma_{\fa,\fm}$.
The definitions of $\tphi_1$ and $\tsigma_1(z)$ mean that
\beql{eq:3N}
\mbox{
$\fa=x\fp_1$ for some $x\in k^\times$ with $x\equiv 1
\bmod^\times \tfm^{(1)}$
}
\eeq
and $\sigma_{\fa,\tfm}$ acts on $\mu(n)$ by $\zeta\mapsto\zeta^z$.
By explicit class field theory over $\bbQ$, this latter condition
translates as
\beql{eq:3O}
N\fa\equiv z \bmod p^{n+1}\bbZ_p
\eeq
We need to show that
$
\{\xi_{\rho,\us,\um}|_{\fa\fp^\um J},\fa\fp^\um J\}_\fm=
\{\xi'_{\rho,\us',\um'},\fp^{\um'} J\}_\fm
$
where $\us':=(zs_1,s_2,\ldots,s_d)$ and
$\um'=(1+m_1,m_2,\ldots,m_d)$ so that
$\fp^{\um'}J=\fp_1\fp^\um J$.
Now $x\fp^{\um'}J=\fa\fp^\um J\subset \fp^\um J$
so we may define $\xi':\fp^{\um'}J\rightarrow \bbC^\times$ by
\[
\xi'(a)=\xi_{\rho,\us,\um}(xa) \ \mbox{for all $a\in\fp^{\um'}J$}
\]
But\refeq{3N} implies $x\in k_{\uinf}^\times$, so\refeq{3M.5} gives
$
\{\xi_{\rho,\us,\um}|_{\fa\fp^\um J},\fa\fp^\um J\}_\fm=
\{\xi',\fp^{\um'} J\}_\fm
$.
It therefore suffices to show that $\xi'=\xi_{\rho,\us',\um'}$, \ie\ that
\beql{eq:3P}
\rho(xa)\zeta_n^{\sum_{i=1}^d s_ip^{-m_i}j\tau_i(x)j\tau_i(a)}=
\rho(a)\zeta_n^{zs_1p^{-1-m_1}j\tau_1(a)+\sum_{i=2}^d s_ip^{-m_i}j\tau_i(a)}
\ \ \ \forall\,a\in\fp^{\um'}J
\eeq
Now $px\in\fa\subset\cO$ and\refeq{3N} implies $px\equiv p\bmod \ff'$
so that $\rho(xa)^p=\rho(pxa)=\rho(pa)=\rho(a)^p$ in $\mu_{f'}$,
hence
\beql{eq:3Q}
\rho(xa)=\rho(a)
\eeq
equation\refeq{3N} also says that if $i\geq 2$ then
$\ord_{\fp_i}(x-1)\geq n+1$ so that
\beql{eq:3R}
j\tau_i(x)\in 1+p^{n+1}\bbZ_p\ \forall\,i\geq 2
\eeq
and since $x\in k_{\uinf}^\times$, equation\refeq{3O}
gives
\beql{eq:3S}
z\equiv p|N_{k/\bbQ}x|= pN_{k/\bbQ}x=
pj\tau_1(x)j\tau_2(x)\ldots j\tau_d(x)
\equiv pj\tau_1(x) \bmod p^{n+1}\bbZ_p
\eeq
equations\refeq{3Q},\refeq{3R} and\refeq{3S} together establish\refeq{3P}, as required.\ePf
\noindent Let $\bbC_p[[\uX]]:=\bbC_p[[X_1,\ldots,X_d]]$
denote the ring of $d$-variable formal power series over
$\bbC_p$. A typical element is written as
$F=\sum_{\um\geq 0}a_\um \uX^\um:=
\sum_{\um\geq 0}a_\um X_1^{m_1}\ldots X_d^{m_d}$
where $\um$ runs through $\bbZ_{\geq 0}^d$. For
any $c$ in $\bbR_{\geq 0}$ we define the  norm $||\cdot||_c$ by setting
\[
||F||_c:=\sup\,\{|a_\um|_pc^{m_1+\ldots+m_d}\,:\,\um\geq 0\}
\in\bbR_{\geq 0}\cup\{\infty\}
\]
Then the set
$\cB(n):=\{F\in\bbC_p[[\uX]]\,:\,\mbox{$||F||_c<\infty$ for some
$c>p^{-1/(p-1)p^n}$ }\}$ is precisely the $\bbC_p$-algebra of those
power series converging and bounded on some `open' polydisc
about $0\in\bbC_p^d$ containing $\mu(n)^d$. The latter group acts on
$\cB(n)$ by setting
\[
\uzeta\cdot F(\uX):=F(\zeta^{(1)}(X_1+1)-1,\ldots,\zeta^{(d)}(X_d+1)-1)\ \ \
\forall\,\uzeta=(\zeta^{(1)},\ldots,\zeta^{(d)})\in\mu(n)^d, F\in\cB(n)
\]
(A rigorous
treatment and discussion of the case $d=2$ is given
in~\cite[Sec. 3.2]{p-adic lims} and the case of general $d$
may be treated in the same way, \emph{mutatis mutandis}). For any
$\ur=(r_1,\ldots,r_d)\in\bbZ_{\geq 0}^d$ such that
$r_i\leq n+1\ \forall\,i$ we write
$\mu(\ur-\uo)$ for the subgroup $\mu(r_1-1)\times\ldots\times\mu(r_d-1)$
of $\mu(n)^d$ and define a
$\bbC_p$-linear projection operator $V_\ur$ onto the subspace of $\cB(n)$ fixed by
$\mu(\ur-\uo)$ by acting with the corresponding idempotent:
\displaymapdef{V_\ur}{\cB(n)}{\cB(n)}{F(\uX)}
{\displaystyle\frac{1}{p^{r_1+\ldots+r_d}}\sum_{\uzeta\in
\mu(\ur-\uo)}\uzeta\cdot F(\uX)}
(The notation has been chosen in order that $V_{\underline{0}}$ should be the
identity).
For any $\uz=(z_1,\ldots,z_d)\in\bbQ_p^d$ we shall use the
abbreviation $(1+\uX)^{\uz}$ for the
product of the binomial series $(1+X_1)^{z_1}\ldots(1+X_d)^{z_d}\in \bbQ_p[[\uX]]$
which lies in $\bbZ_p[[\uX]]$ (hence in $\cB(n)$) if $\uz\in\bbZ_p^d$.
\begin{lemma}\label{lemma:3H}\ \vspace{1ex}\\
(i)~If $\uz\in\bbZ_p$ then
\[
V_\ur (1+\uX)^{\uz}=
\left\{
\begin{array}{ll}
(1+\uX)^{\uz}&
\mbox{if $p^{r_i}|z_i$ for all $i$ and}\\
0&
\mbox{otherwise}
\end{array}
\right.
\]
(ii)~If $M$ is any (topologically) closed subfield of $\bbC_p$ with
ring of valuation integers
$\cO_M$ then $V_\ur\cO_M[[\uX]]\subset\cO_M[[\uX]]$.
\end{lemma}
\bPf\ The group $\mu(n)^d$ acts on
$(1+\uX)^{\uz}$ by multiplication by the character sending
$\uzeta$ to ${\zeta^{(1)}}^{z_1}\ldots{\zeta^{(d)}}^{z_d}$. This character
is trivial on $\mu(\ur-\uo)$ if and only if $p^{r_i}|z_i$ for all $i$.
Part~(i) follows as does the statement of part~(ii) with $\cO_N[[\uX]]$
replaced by $\cO_N[\uX]$ (since the latter is clearly
spanned over
$\cO_N$ by $\{(1+\uX)^{\um}\,:\,\um\in\bbZ_{\geq 0}^d\}$).
Now if $F$ lies in $\cO_N[[\uX]]$ then $||F||_c<\infty$
for any $c$ with
$1>c>p^{-1/(p-1)p^n}$. Moreover, for such $c$ we can find a sequence of
polynomials
$F_l\in\cO_N[\uX]$ tending to $F$ w.r.t. $||\cdot||_c$
as $l\rightarrow\infty$ (\eg\ a sequence of truncations of $F$).
On the other hand, the $\mu(n)^d$ action preserves the norm $||\cdot||_c$
(see \eg\ \cite[Sec. 3.2]{p-adic lims}) from which
it follows that $V_\ur$ is $||\cdot||_c$-continuous.
Thus $V_\ur F_l$ tends to $V_\ur F$ with respect to $||\cdot||_c$,
hence also coefficientwise. But $V_\ur F_l$ lies in
$\cO_M[\uX]\ \forall\,l$ so
$V_\ur F$ lies in $\cO_M[[\uX]]$ .\ePf
\noindent The following Proposition is crucial to the proof of
Claim~\ref{claim:3B} and of some interest in its own right. For ease of
reading we shall defer its proof to the end of this section.
\begin{prop}\label{prop:3B}
Suppose given $\rho: J\rightarrow\bbC^\times$ with
${\rm ann}_\cO(\rho)=\ff'$ as in Lemma~\ref{lemma:3G}. Then
there exists a family $\{F_{\um}\,:\,\um\in\bbZ^d\}$ of
$d$-variable, formal $p$-adic power series (depending on $\rho$ and
$j$) such that
\begin{enumerate}
\item\label{A'} $F_{\um}(\uX)$ lies in $\bbZ_p[\mu_{f'}][[\uX]]$
for all $\um\in\bbZ^d$,
\item\label{C'} $V_\ur F_{\um}(\uX)=
F_{\ur+\um}((1+X_1)^{p^{r_1}}-1,\ldots,(1+X_d)^{p^{r_d}}-1)$
for all $\um\in\bbZ^d$ and $\ur=(r_1,\ldots,r_d)\in\bbZ_{\geq 0}^d$ with
$r_i\leq n+1\ \forall\,i$, and
\item\label{B'} $j(Z(0;\fw_{\rho,\us,\um}))=
F_{\um}(\zeta_n^{s_1}-1,\ldots,\zeta_n^{s_d}-1)$
for all $\um\in\bbZ^d$ and $\us\in\cS(\un)$.
\end{enumerate}
\end{prop}
Returning now to equation\refeq{3M.3} we choose a fractional
ideal $J$ of $k$ and a character of finite order $\xi:J\rightarrow\bbC^\times$
with ${\rm ann}_\cO(\xi)=\ff$ and such that
$h\cdot\fw_\fm^0=\{\xi,J\}_\fm$ in $\fW_\fm$. Weak approximation
and the equivalence
relation\refeq{3M.5} allow us to assume
that $J$ is prime to $p$ so that, by Lemma~\ref{lemma:3F}, we have
$\xi=\xi_{\rho,\us,\uzero}$ for some
$\rho: J\rightarrow\bbC^\times$ with ${\rm ann}_\cO(\rho)=\ff'$
and $\us\in\cS(\un)$. Thus $h\cdot\fw_\fm^0=\fw_{\rho,\us,\uzero}$.
Let $\{F_{\um}\,:\,\um\in\bbZ^d\}$ be a family of power
series satisfying in Proposition~\ref{prop:3B} w.r.t.\ $J$ and $\rho$.
Using also Lemma~\ref{lemma:3G} \etc, equation\refeq{3M.3} becomes
\begin{eqnarray*}
a_{h,\underline{t}}&=& p^{-t_1-\ldots-t_d}
\sum_{z_i\in Z_i\atop i=1,\ldots,d}
\sum_{l_i\in L_i\atop i=1,\ldots,d}
\left(
j(Z(0;\fw_{\rho,\uz\us,\ul-\ut}))\prod_{i=1}^d
\hat{\sigma}_i(z_i)\hat{\phi}^{l_i}
(g_i(\zeta_{n_i-t_i}-1))
\right)\\
&=& p^{-t_1-\ldots-t_d}
\sum_{l_i\in L_i\atop i=1,\ldots,d}
\sum_{z_i\in Z_i\atop i=1,\ldots,d}
\left(
F_{\ul-\ut}(\zeta_n^{z_1s_1}-1,\ldots,\zeta_n^{z_ds_d}-1)
\prod_{i=1}^d
\hat{\phi}^{l_i}
(g_i(\zeta_{n_i-t_i}^{z_i}-1))
\right)
\end{eqnarray*}
Since $t_i$ is strictly less than $n_i+1$ it follows that, for each $i$,
$\mu(t_i-1)$ acts (multiplicatively and freely) on the difference set
$\mu(n_i)\setminus\mu(n_i-1)$ of
\emph{primitive} $p^{n_i+1}$th roots of unity
and we choose a set $\Gamma(n_i,t_i)$ of orbit representatives
for this action for each $i$. Fixing $i$ we see that
$\zeta_{n_i}^{z_i}$ lies in $\mu(n_i)\setminus\mu(n_i-1)$
for every $z_i\in Z_i$ and so may be written
\[
\zeta_{n_i}^{z_i}=\gamma_i\nu_i\ \ \ \mbox{for some unique
$\gamma_i\in \Gamma(n_i,t_i)$ and $\nu_i\in\mu(t_i-1)$ (depending on $z_i$)}
\]
Thus $\zeta_n^{z_is_i}=\zeta_{n_i}^{z_is'_i}=
\gamma_i^{s'_i}\nu_i^{s'_i}$ where we write $s'_i$ for
$p^{n_i-n}s_i\in\bbZ_p^\times$.
Moreover, as $z_i$ runs once through $Z_i$, so
$\zeta_{n_i}^{z_i}$ runs once through  $\mu(n_i)\setminus\mu(n_i-1)$,
hence the pair
$(\gamma_i,\nu_i)$ runs once through $\Gamma(n_i,t_i)\times \mu(t_i-1)$,
and the same must obviously be true of the pair $(\gamma_i,\nu_i^{s'_i})$.
On the other hand $\zeta_{n_i-t_i}^{z_i}=
(\zeta_{n_i}^{z_i})^{p^{t_i}}=\gamma_i^{p^{t_i}}$
which is independent of $\nu_i$ and clearly runs exactly once through
$\mu(n_i-t_i)\setminus\mu(n_i-t_i-1)$
as $\gamma_i$ runs
through $\Gamma(n_i,t_i)$. Putting this together
and using Proposition~\ref{prop:3B}~\ref{C'},
the last equation gives
\begin{eqnarray*}
a_{h,\underline{t}}&=&\sum_{l_i\in L_i\atop i=1,\ldots,d}
\sum_{\gamma_i\in \Gamma(n_i,t_i)\atop i=1,\ldots,d}
\sum_{\nu_i\in\mu(t_i-1)\atop i=1,\ldots,d}
\left(
p^{-t_1-\ldots-t_d}
F_{\ul-\ut}(\gamma_1^{s'_1}\nu_1^{s'_1}-1,\ldots,
\gamma_d^{s'_d}\nu_d^{s'_d}-1)
\prod_{i=1}^d\hat{\phi}^{l_i}(g_i(\gamma_i^{p^{t_i}}-1))
\right)
\\
                   &=&
\sum_{l_i\in L_i\atop i=1,\ldots,d}
\sum_{\gamma_i\in \Gamma(n_i,t_i)\atop i=1,\ldots,d}
\left(
V_\ut
F_{\ul-\ut}(\gamma_1^{s'_1}-1,\ldots,
\gamma_d^{s'_d}-1)
\prod_{i=1}^d\hat{\phi}^{l_i}(g_i(\gamma_i^{p^{t_i}}-1))
\right)
\\
                   &=&
\sum_{l_i\in L_i\atop i=1,\ldots,d}
\sum_{\gamma_i\in \Gamma(n_i,t_i)\atop i=1,\ldots,d}
\left(
F_{\ul}(\gamma_1^{s'_1p^{t_1}}-1,\ldots,
\gamma_d^{s'_dp^{t_d}}-1)
\prod_{i=1}^d\hat{\phi}^{l_i}(g_i(\gamma_i^{p^{t_i}}-1))
\right)
\\
                   &=&
\sum_{l_i\in L_i\atop i=1,\ldots,d}
\sum_{\zeta^{(i)}\in\mu(n_i-t_i)\setminus\mu(n_i-t_i-1)\atop i=1,\ldots,d}
\left(
F_{\ul}( \zeta^{(1)s'_1}-1,\ldots,
\zeta^{(d)s'_d}-1)
\prod_{i=1}^d\hat{\phi}^{l_i}(g_i(\zeta^{(i)}-1))
\right)
\end{eqnarray*}
Substituting  this equation into\refeq{3M} and noting that $g_i(0)=0$ (by
Lemma~\ref{lemma:3E}~\ref{B}) we get
\begin{eqnarray}
a_h&=&p^{-N}
\sum_{\zeta^{(i)}\in\mu(n_i)\setminus\mu(-1)\atop i=1,\ldots,d}
\sum_{l_i\in L_i\atop i=1,\ldots,d}
\left(
F_{\ul}( \zeta^{(1)s'_1}-1,\ldots,
\zeta^{(d)s'_d}-1)
\prod_{i=1}^d\hat{\phi}^{l_i}(g_i(\zeta^{(i)}-1))
\right)\nonumber\\
  &=&p^{-N}\sum_{\zeta^{(i)}\in\mu(n_i)\atop i=1,\ldots,d}
\sum_{l_i\in L_i\atop i=1,\ldots,d}
\left(
F_{\ul}( \zeta^{(1)s'_1}-1,\ldots,
\zeta^{(d)s'_d}-1)
\prod_{i=1}^d\hat{\phi}^{l_i}(g_i(\zeta^{(i)}-1))
\right)\nonumber\\
  &=&
V_{\un+\uo}F(0,\ldots,0)\label{eq:3S.2}
\end{eqnarray}
where we have set $\un+\uo:=(n_1+1,\ldots,n_d+1)$ and
\beql{eq:3S.3}
F(\uX):=\sum_{l_i\in L_i\atop i=1,\ldots,d}
\left(
F_{\ul}((1+X_1)^{s'_1}-1,\ldots,(1+X_d)^{s'_d}-1 )
\prod_{i=1}^d\hat{\phi}^{l_i}\cdot g_i(X_i)
\right)
\eeq
The power series $(1+X_i)^{s'_i}-1$ has coefficients in
$\bbZ_p$ for all $i$ so it follows from Lemma~\ref{lemma:3E}~\ref{A}
and Proposition~\ref{prop:3B}~\ref{A'} that $F(\uX)$ lies in
$\cO_H[\mu_{f'}][[\uX]]$. (As well as $J$,$\rho$ and $\underline{s}$, the
power series $F$ depends on the choices of the $L_i$, the
$F_\ul(\uX)$ for $\ul\in L_1\times\ldots\times L_d$ and the $g_i$.)
Consequently $a_h$ lies in $\cO_H[\mu_{f'}]\subset\bar{\bbZ}_p$
(hence actually in $\bbZ_p$) by
Lemma~\ref{lemma:3H}~(ii).
Thus Claim~\ref{claim:3B} will be established once we have
proved Proposition~\ref{prop:3B}.

We shall shortly see that the existence of the power series $F_{\um}(\uX)$ appearing in
Proposition~\ref{prop:3B} is a  very
natural consequence of the method of
Shintani~\cite{Shi} (with refinements from~\cite{Colm})
for evaluating Dirichlet series like $Z(s;\fw_{\rho,\us,\um})$ by means
of cone decompositions of fundamental domains
for the action of units on $k^\times_{\uinf}$.
Such decompositions are best visualised by using the embedding
$\utau: k\rightarrow \bbR^d$ where $\utau(a)$ is defined to be
$(\tau_1(a),\ldots,\tau_d(a))$. This
extends to an isomorphism $k\otimes_\bbQ \bbR\cong\bbR^d$ and
sends $\utau(k^\times_{\uinf})$ into $\bbR_{>0}^d$.
Similarly, we define an embedding $\ujct: k\rightarrow \bbQ_p^d$ by
$\ujct(a):=(j\tau_1(a),\ldots,j\tau_d(a))$
(which extends to an isomorphism $k\otimes_\bbQ\bbQ_p\cong \bbQ_p^d$).
For any $\um\in\bbZ^d$ we shall write $p^{-\um}\ujct(a)$ for the $d$-tuple
$(p^{-m_1}j\tau_1(a),\ldots,p^{-m_d}j\tau_d(a))$. Recall also the notation
$\fp^\um$ for $\fp_1^{m_1}\ldots\fp_d^{m_d}$
and the fact that the character $\rho$ appearing in Proposition~\ref{prop:3B} has been uniquely
extended to a homomorphism from
$\bbZ[{\textstyle\frac{1}{p}}]J=
\bigcup_{\um\in\bbZ^d}\fp^\um J$ onto $\mu_{f'}\subset\bbC^\times$.
We shall write simply $j\rho$ for the composite $j\circ\rho$.
Since we are assuming that
$J$ is prime to $p$ and $j\tau_i$ induces $\fp_i$ for each $i$,
it follows that
\beql{eq:3S.5}
\mbox{If $a\in\bbZ[{\textstyle\frac{1}{p}}]J$
and $\um\in\bbZ^d$ then
$a\in\fp^\um J\Leftrightarrow\ujct(a)\in
{\displaystyle \prod_{i=1}^d}\left(p^{m_i}\bbZ_p\right)
\Leftrightarrow p^{-\um}\ujct(a)\in \bbZ_p^d$
}
\eeq
Now, for any $d'$-tuple  $\bv=(v_1,\ldots,v_{d'})$ of $\bbQ$-linearly independent
elements of $k^\times_{\uinf}$ (with $1\leq d'\leq d$) we
write  $C(\bv)$ for the open simplicial cone in $\bbR_{>0}^d$ spanned by
$\utau(v_1),\ldots,\utau(v_{d'})$, namely $C(\bv)=\sum_{l=1}^{d'}\bbR_{>0}\utau(v_l)$.
Given $\rho$, we make the following informal \emph{pseudo-definition}:
\beql{eq:3T}
\mbox{``
$
F_{\um,\bv}(\uX)={\displaystyle \sum_{a\in \fp^{\um} J\atop \utau(a)\in C(\bv)}}
j\rho(a)
(1+\uX)^{p^{-\um}\ujct(a)}
$
'' for every $\um\in \bbZ^d$}
\eeq
The R.H.S. of\refeq{3T} is an infinite sum of formal power series in $\bbC_p[[\uX]]$
each of whose constant terms is a root of unity.
Since there is no obvious sense in which it converges,
the value of equation\refeq{3T} is largely conceptual.
To give a meaningful definition of $F_{\um,\bv}(\uX)$ we note first that,
for any $\um\in\bbZ^d$ and $l\in\{1,\ldots,d'\}$,
the intersection $\bbQ_{>0} v_l\cap\fp^\um J$
equals $\bbZ_{>0}w^0_{\um,l}$
for a unique element $w^0_{\um,l}$ of $\ktuinf$.
Now choose \emph{any} element $\bw=(w_1,\ldots,w_{d'})$ of
$(\ktuinf)^{d'}$ such that $w_l\in\bbQ_{>0} v_l\cap\fp^\um J$
(for instance $\bw=\bw^0_\um:=(w^0_{\um,1},\ldots,w^0_{\um,d'})$).
Let us denote by $P(\bw)$ the half-open parallelipiped
$\sum_{l=1}^{d'}(0,1]\utau(w_l)\subset C(\bv)$
and define $F_{\um,\bv}(\uX)$ to be $G_{\um,\bw}(\uX)/H_{\um,\bw}(\uX)$ where
\[
G_{\um,\bw}(\uX):=
{\displaystyle \sum_{a\in\fp^{\um} J\atop \utau(a)\in P(\bw)}}
j\rho(a)(1+\uX)^{p^{-\um}\ujct(a)}
\]
and
\[
H_{\um,\bw}(\uX):=
{\displaystyle \prod_{l=1}^{d'}}
\left(
1-j\rho(w_l)(1+\uX)^{p^{-\um}\ujct(w_l)}
\right)
\]
The sum defining $G_{\um,\bw}(\uX)$ is finite since
$P(\bw)$ has compact closure and
$\utau(\fp^{\um} J)$ is discrete in $\bbR^d$.
Since $H_{\um,\bw}(\uX)$ is non-zero, we may certainly consider $F_{\um,\bv}(\uX)$
as an element of the fraction field
of $\bbC_p[[\uX]]$ which might, \emph{a priori}, depend on the choice of $\bw$.
In fact:
\begin{lemma}\label{lemma:3I} For $\rho$, $\bv$, $\um$ and $\bw$ as above
\begin{enumerate}
\item\label{part:lemma3I1} $G_{\um,\bw}(\uX)$ and $H_{\um,\bw}(\uX)$
lie in $\bbZ_p[\mu_{f'}][[\uX]]$.
\item\label{part:lemma3I2} $F_{\um,\bv}(\uX)$ is independent of
the choice of $\bw$ (given $\rho$, $\bv$ and $\um$).
\item\label{part:lemma3I3} If $\bbQ v_l\cap \fp^\um J\not\subset\ker(\rho)$
for $l=1,\ldots,d'$ then $F_{\um,\bv}(\uX)$ lies in $\bbZ_p[\mu_{f'}][[\uX]]$.
\end{enumerate}
\end{lemma}
\bPf Part~\ref{part:lemma3I1} follows from equation\refeq{3S.5}.
Let us write $w_l$ as $q_lw^0_{\um,l}$ with $q_l\in \bbZ_{>0}$ for $l=1,\ldots,d'$
and set $R(\uX)=
\prod_{l=1}^{d'}
\left(\sum_{r_l=0}^{q_l-1}
j\rho(r_lw^0_{\um,l})(1+\uX)^{p^{-\um}\ujct(r_lw^0_{\um,l})}\right)$.
Then clearly, $H_{\um,\bw}(\uX)=H_{\um,\bw^0_{\um}}(\uX)R(\uX)$. Furthermore,
$G_{\um,\bw}(\uX)=G_{\um,\bw^0_{\um}}(\uX)R(\uX)$ because
$w^0_{\um,l}$ lies in $\fp^{\um} J\ \forall l$ and so the set
$\{a\in\fp^{\um} J\,:\,\utau(a)\in P(\bw)\}$
is the disjoint union of the translates
$\{a\in\fp^{\um} J\,:\,\utau(a)\in P(\bw^0_{\um})\}+
(r_1w^0_{\um,1}+\ldots+r_{d'}w^0_{\um,d'})$
as $r_l$ runs from $0$ to $q_l-1$ for all $l$.
Thus $G_{\um,\bw}(\uX)/H_{\um,\bw}(\uX)=G_{\um,\bw^0_{\um}}(\uX)/H_{\um,\bw^0_{\um}}(\uX)$
which proves~\ref{part:lemma3I2}.
Finally, if  $\bbQ v_l\cap \fp^\um J\not\subset\ker(\rho)\ \forall l$
then for each $l$ we may take $w_l\in\bbQ_{>0} v_l\cap \fp^\um J$ such that
$j\rho(w_l)$ is a
\emph{non-trivial} $f'$-th root of unity in $\bbC_p$.
Since $p\ndiv f'$, the constant term
$\prod_{l=1}^{d'}(1-j\rho(w_l))$ of $H_{\um,\bw}(\uX)$
is then a \emph{unit} of $\bbZ_p[\mu_{f'}]$ so that
$H_{\um,\bw}(\uX)$ is a unit of $\bbZ_p[\mu_{f'}][[\uX]]$ and
part~\ref{part:lemma3I3} follows from part~\ref{part:lemma3I1}.\ePf
\rem\ The pseudo-formula\refeq{3T} can be seen as the result of an informal limiting process
in which the integers $q_l$ in the above proof tend to infinity in such a
way that $H_{\um,\bw}(\uX)$ tends to $1$.\vspace{1ex}\\
In view of Lemma~\ref{lemma:3I}~\ref{part:lemma3I3}
we introduce the following hypothesis which may or may not be verified by
a character $\rho:J\rightarrow\bbC^\times$
and a $d'$-tuple $\bv$ as above:
\vspace{1ex}\\
{\bf Hypothesis}
For each $l=1,\ldots,d'$ we have
H($\rho$,$\bv$): We have $\bbQ v_l\cap J\not\subset\ker(\rho)$.
\vspace{1ex}\\
If $\um\in\bbZ^d$ then
$p^N(\bbQ v_l\cap J)\subset \bbQ v_l\cap\fp^\um J\subset p^{-N}(\bbQ v_l\cap J)$ for
some $N>>0$ and since $\rho$ takes values in $\mu_{f'}$ and $p\ndiv f'$,
it follows that \Hrhov\ is equivalent
to $\bbQ v_l\cap \fp^\um J\not\subset\ker(\rho)$
for all $l=1,\ldots,d'$ \emph{and all} $\um\in\bbZ^d$
(which in turn is equivalent to
$\rho(w^0_{l,\um})\neq 1$ for all $l$ and $\um$).
\begin{lemma}\label{lemma:3I.5}
Suppose \Hrhov\ holds. Then
$F_{\um,\bv}(\uX)\in\bbZ_p[\mu_{f'}][[\uX]]\ \forall\,\um\in\bbZ^d$.
Moreover,
for any $\ur=(r_1,\ldots,r_d)\in\bbZ_{\geq 0}^d$ with
$r_i\leq n+1\ \forall\,i$, we have
$V_\ur F_{\um,\bv}(\uX)=
F_{\ur+\um,\bv}((1+X_1)^{p^{r_1}}-1,\ldots,(1+X_d)^{p^{r_d}}-1)
\ \forall\,\um\in\bbZ^d$.
\end{lemma}
\bPf\ The first statement follows from
Lemma~\ref{lemma:3I}~\ref{part:lemma3I3} and the above discussion. For the
second, choose $\bw$ such that $w_l$ lies in $\bbQ_{>0}v_l\cap
\fp^{\ur+\um}J$ for all $l=1,\ldots,d'$. Then
$G_{\um,\bw},H_{\um,\bw},G_{\ur+\um,\bw}$ and $H_{\ur+\um,\bw}$
are all defined and lie in $\bbZ_p[\mu_{f'}][[\uX]]$. Moreover,
it follows from\refeq{3S.5} (with $\ur+\um$ in place of $\um$)
that $p^{-m_i}j\tau_i(w_l)\in p^{r_i}\bbZ_p$ for all $l$ and
$i$, with the result that $\mu(\ur-\uo)$ acts trivially on $H_{\um,\bw}(\uX)$.
Therefore
\beql{eq:3V} H_{\um,\bw}(\uX)V_\ur
F_{\um,\bv}(\uX)=V_\ur(H_{\um,\bw}(\uX)F_{\um,\bv}(\uX))= V_\ur
G_{\um,\bw}(\uX) \eeq
If $a\in\fp^{\um}J$ then
Lemma~\ref{lemma:3H}~(i) combined with\refeq{3S.5} shows that
$V_\ur (1+\uX)^{p^{-\um}\ujct(a)}$ equals
$(1+\uX)^{p^{-\um}\ujct(a)}$ if $a\in\fp^{\ur+\um}$ and $0$ otherwise. It
follows easily from the definitions of $G_{\um,\bw}(\uX)$ and
$H_{\um,\bw}(\uX)$ that
\begin{eqnarray}
V_\ur G_{\um,\bw}(\uX)&=&G_{\ur+\um,\bw}
                             ((1+X_1)^{p^{r_1}}-1,\ldots,(1+X_d)^{p^{r_d}}-1)\nonumber\\
          &=&H_{\ur+\um,\bw}((1+X_1)^{p^{r_1}}-1,\ldots,(1+X_d)^{p^{r_d}}-1)\times\nonumber\\
          &&\ \ \ F_{\ur+\um,\bv}((1+X_1)^{p^{r_1}}-1,\ldots,(1+X_d)^{p^{r_d}}-1)\nonumber\\
          &=&H_{\um,\bw}(\uX)
                F_{\ur+\um,\bv}((1+X_1)^{p^{r_1}}-1,\ldots,(1+X_d)^{p^{r_d}}-1)\nonumber
\end{eqnarray}
Combining the last equation with\refeq{3V} and cancelling
$H_{\um,\bw}(\uX)$ gives the result.\ePf
\rem\ The inequalities $r_i\leq n+1$ in the statement of the Lemma
are unnatural and unnecessary. Indeed, if $n$ is
allowed to tend to infinity, we get an action of $\mu(\infty)$ on the subring
$\bigcap_{n\geq 0}\cB(n)$ of $\bbC_p[[X]]$ (namely, those
power series converging on the `open' polydisc radius $1$ about $0\in\bbC_p^d$) which contains
$\bbZ_p[\mu_{f'}][[\uX]]$. The argument then goes through for any $\ur\in\bbZ_{\geq 0}^d$.
Furthermore, the second statement of the Lemma is easily seen to hold even
when \Hrhov\ fails, provided we use a natural extension
of the $\mu(\infty)$ action to the fraction
field of $\bbZ_p[\mu_{f'}][[\uX]]$.\vspace{1ex}\\
Recall that for any
$\rho:J\rightarrow \bbC^\times$, $\us\in\cS(n)$ and $\um\in\bbZ^d$
as in Propostion~\ref{prop:3B} we have defined a character
$\xi_{\rho,\us,\um}:\fp^\um J\rightarrow\bbC^\times$ by\refeq{3M.75}.
This in turn defines a complex Dirichlet series $Z_{\bv}(s,\xi_{\rho,\us,\um},\fp^\um J)$
by the following formula.
\[
Z_{\bv}(s,\xi_{\rho,\us,\um},\fp^\um J):=
{\displaystyle \sum_{a\in \fp^\um J\atop \utau(a)\in C(\bv)}}
\frac{\xi_{\rho,\us,\um}(a)}{(\tau_1(a)\ldots\tau_d(a))^s}=
{\displaystyle \sum_{a\in \fp^\um J\atop \utau(a)\in C(\bv)}}
\xi_{\rho,\us,\um}(a) N_{k/\bbQ}(a)^{-s}
\]
\begin{lemma}\label{lemma:3J} Suppose \Hrhov\ holds for some $d'$-tuple $\bv$ as above.
Then for any $\um\in\bbZ^d$ and $\us\in\cS(n)$
the function
$Z_{\bv}(s,\xi,I)$ converges absolutely for ${\rm Re}(s)>d'/d$
and possesses a meromorphic continuation to $\bbC$.
Moreover the latter is regular at zero with value
in $\bbQ(\mu_{f'p^{n+1}})$. More precisely, we have
\beql{eq:3V.5}
j(Z_{\bv}(0,\xi_{\rho,\us,\um},\fp^\um J))=F_{\um,\bv}(\zeta_n^{s_1}-1,\ldots,\zeta_n^{s_d}-1)
\eeq
\end{lemma}
\bPf\ This all follows from the results of~\cite{Shi},
the particularly simple form of\refeq{3V.5} being a consequence our Hypothesis \Hrhov.
We sketch the derivation.
Let us fix $\um$, $\us$ and
choose $\bw=(w_1,\ldots,w_{d'})$ with
$w_l\in\bbQ_{>0} v_l\cap \fp^\um J$ for $l=1,\ldots,d'$ and $\rho(w_l)\neq 1\ \forall\,l$
(by \Hrhov).
The set
$\{a\in\fp^{\um} J\,:\,\utau(a)\in C(\bv)\}$
is the disjoint union of the translates
$\{a\in\fp^{\um} J\,:\,\utau(a)\in P(\bw)\}+
(r_1w_1+\ldots+r_{d'}w_{d'})$
as $r_l$ runs from $0$ to $\infty$ for each $l$.
It follows that
$Z_{\bv}(s,\xi_{\rho,\us,\um},\fp^\um J)=
\sum_a \xi_{\rho,\us,\um}(a)Z_{\bw}(s,a,\xi_{\rho,\us,\um})$
where $a$ runs through the finite set
$\{a\in\fp^{\um} J\,:\,\utau(a)\in P(\bw)\}$
and
\[
Z_{\bw}(s,a,\xi_{\rho,\us,\um}):=\sum_{r_1,\ldots r_{d'}\in\bbZ_{\geq 0}}
\frac{\xi_{\rho,\us,\um}(w_1)^{r_1}\ldots \xi_{\rho,\us,\um}(w_{d'})^{r_{d'}}}
{\prod_{i=1}^d\tau_i(a+(r_1w_1+\ldots+r_{d'}w_{d'}))^s}
\]
The last two equations are formal until the absolute convergence of
$Z_{\bw}(s,a,\xi_{\rho,\us,\um})$ is established. But this series is just
`$\zeta(s,A,x,\chi)$' in the notation of~\cite[p.~396]{Shi} where Shintani's
$(n,r)$ is our $(d,d')$, his $A$ is the matrix
$(a_{li}):=(\tau_i(w_l))$ with $1\leq l\leq d'$, $1\leq i\leq d$, his $x$ is the $d'$-tuple
$(x_1,\ldots,x_{d'})\in (\bbQ\cap (0,1])^{d'}$ such that $a=x_1w_1+\ldots+x_{d'}w_{d'}$
and his $\chi$ is the $d'$-tuple
$(\chi_1,\ldots,\chi_{d'}):=(\xi_{\rho,\us,\um}(w_1),\ldots,\xi_{\rho,\us,\um}(w_{d'}))$.
The first statement of the lemma now follows from that of~\cite{Shi}, Proposition~1.
Furthermore, this Proposition actually implies the regularity of
$Z_{\bw}(s,a,\xi_{\rho,\us,\um})$
at $s={1-m}$ for all
$m\in\bbZ_{\geq 1}$ and (taking $m=1$) the formula
\[
Z_{\bw}(0,a,\xi_{\rho,\us,\um})=\frac{1}{d}\sum_{i=1}^d B_1(A,1-x,\chi)^{(i)}
\]
where $B_1(A,1-x,\chi)^{(i)}$ denotes the constant term of a certain Laurent series.
But since $\rho(w_l)$ is a \emph{non-trivial} $f'$th root of unity and $p\ndiv f'$,
it follows from\refeq{3M.75} that $\chi_l\neq 1\ \forall\,l$
which means that these Laurent series are actually all power series
in Shintani's variables $u$ and $t_i,\ldots,t_{i-1},t_{i+1},\ldots,t_d$. Setting these
all equal to zero we deduce that $B_1(A,1-x,\chi)^{(i)}$ is simply equal to
$(\prod_{l=1}^{d'}(1-\chi_l))\inv$ independently of $i$ and $a$,
which is therefore also the value of $Z_{\bw}(0,a,\xi_{\rho,\us,\um})$ for all $a$. Thus
\[
Z_{\bv}(0,\xi_{\rho,\us,\um},\fp^\um J)=
\frac
{ \sum_{a\in\fp^{\um} J\atop \utau(a)\in P(\bw)}
\xi_{\rho,\us,\um}(a)}
{ \prod_{l=1}^{d'}
(1-\xi_{\rho,\us,\um}(w_l))}
\]
The result follows from the definition of $F_{\um,\bv}$ on applying $j$ to both sides
of the above equation, noting that if $a\in\fp^\um J$ then
$j(\xi_{\rho,\us,\um}(a))$ is the value of the power series
$j\rho(a)(1+\uX)^{p^{-\um}\ujct(a)}$ at $\uX=(\zeta_n^{s_1}-1,\ldots,\zeta_n^{s_d}-1)$,
by\refeq{3M.75}.
(Recall that ``$j(\zeta_n)=\zeta_n$'' in our notation.)\ePf
\noindent
The group $E_\uinf=E_\uinf(k)$ of totally positive units of $k$
acts on $\bbR^d_{>0}$ by setting
$\varepsilon (x_1,\ldots,x_d):=(\tau_1(\varepsilon)x_1,\ldots,\tau_d(\varepsilon)x_d)$
for all $\varepsilon\in E_\uinf\mbox{\ and \ }(x_1,\ldots,x_d)\in\bbR^d_{>0}$.
Given any subgroup
$E$ of finite index in $E_\uinf$, a
\emph{$k$-rational cone decomposition
of a fundamental domain for $E$ acting on $\bbR_{>0}^d$}
(called just a \emph{cone decomposition for $E$} for short) is a finite
set $\cC$ of $d'$-tuples $\bv$ of $\bbQ$-linearly independent
elements of $k^\times_{\uinf}$ (with $d'$ depending on
$\bv$ and varying between $1$ and $d$) such that the cones
$\varepsilon C(\bv)=C(\varepsilon \bv)$ for $\bv\in\cC$ and $\varepsilon\in E$
are pairwise disjoint and their union is
$\bbR_{>0}^d$. Shintani proved in~\cite{Shi} that
cone decompositions exist for any such $E$. Their relevance to our
situation comes from a corresponding decomposition of
twisted zeta-functions. Recall that for $\rho:J\rightarrow \bbC^\times$, $\us\in\cS(n)$ and
$\um\in\bbZ^d$ as above,
Lemma~\ref{lemma:3F} implies that ${\rm ann}_\cO(\xi_{\rho,\us,\um})=\ff$
where $\fm=\ff\infty$. In this situation:
\begin{lemma}\label{lemma:3K}
For any subgroup $E$ of finite index in $E_\fm$
(hence in $E_\uinf$) and any cone decomposition $\cC$ for $E$ we have
the following equality (as meromorphic functions of $s\in\bbC$)
\[
Z(s,\fw_{\rho,\us,\um})=
\frac{N(\fp^\um J)^s}{|E_\fm:E|}\sum_{\bv\in\cC}Z_{\bv}(s,\xi_{\rho,\us,\um},\fp^\um J)
\]
\end{lemma}
\bPf\ By meromorphic continuation, it suffices to prove the
equality for ${\rm Re}(s)>1$. Now by definition, $\fw_{\rho,\us,\um}=
\{\xi_{\rho,\us,\um},\fp^\um J\}_{\fm}\in\fW_\fm$. Therefore, the definitions, notations and
equations of~\cite[pp.\ 15,16]{twizas} give
\begin{eqnarray*}
Z(s,\fw_{\rho,\us,\um})=Z_\emptyset(s,\fw_{\rho,\us,\um})
                       &=&\frac{1}{|E_\fm:E|}Z_\emptyset(s,\fw_{\rho,\us,\um},E)\\
                       &=&\frac{N(\fp^\um J)^s}{|E_\fm:E|}
                          \sum_{a\in \cR}
                           \xi_{\rho,\us,\um}(a) N_{k/\bbQ}(a)^{-s}
\end{eqnarray*}
for all $s\in\bbC$, ${\rm Re}(s)>1$, where $\cR$ is any set of orbit
representatives for the (multiplicative)
action of $E$ on $k_{\uinf}^\times\cap \fp^\um J$. We may clearly take
$\cR$ to be the set $\{a\in\fp^\um J\,:\,\utau(a)\in\bigcup_{\bv\in\cC}C(\bv))\}=
\bigcup_{\bv\in\cC}\{a\in\fp^\um J\,:\,\utau(a)\in C(\bv))\}$
(a disjoint union). The Lemma therefore follows by absolute convergence from the definition
of $Z_{\bv}(s,\xi_{\rho,\us,\um},\fp^\um J)$.\ePf
\noindent The final ingredient in proving Proposition~\ref{prop:3B}
concerns the existence of subgroups of $E_\fm$ and cone decompositions
with the properties we require.
\begin{lemma}\label{lemma:3L} Suppose given $\rho: J\rightarrow\bbC^\times$ with
${\rm ann}_\cO(\rho)=\ff'$. Then there exists a subgroup $E$ of finite index
in $E_\fm$ and a cone decomposition $\cC$ for $E$ satisfying
\beql{eq:3X}
p\ndiv[E_\fm:E]
\eeq
and
\beql{eq:3Y}
\mbox{\Hrhov\ holds for all $v\in\cC$}
\eeq
\end{lemma}
\bPf\ We shall use Colmez' construction of a
cone decomposition in~\cite[\S 2]{Colm}, with two small modifications.
We shall harmonise our set-up with his by suppressing $\utau$ and
regarding $k$ as embedded in $\bbR^d$. A direct application of Colmez'
Lemme~2.1 \emph{loc.\ cit.}
(with $F=k$, $n=d$, $V=E_\fm$) would give rise to elements
$\varepsilon_1,\ldots,\varepsilon_{d-1}$ of $E_\fm$
satisfying Colmez' hypothesis $(H)$.
However, to ensure that they generate a subgroup $E$ satisfying
Condition\refeq{3X}, we need to modify their construction
in the proof of Lemme~2.1 to which we now refer.
Specifically, we recall that Colmez chooses a positive real number
$r(V)=r(E_\fm)$ such that ${\rm Log}(E_\fm)$ has non-empty intersection with the ball
$B(l_t(M),r(E_\fm))$  for each $i$ and each $M>0$.
(We follow him in using the $\sup$-norm $||\cdot||$ on $\bbR^d$.)
Let us choose a  $\bbZ$-basis $\eta_1,\ldots,\eta_{d-1}$ of $E_\fm$
and set
$r'(E_\fm)=r(E_\fm)+(p-1)(d-1){\rm max}\{||{\rm Log}(\eta_t)||\,:\,1\leq t\leq d-1\}$.
It is easy to see that the larger ball $B(l_t(M),r'(E_\fm))$ must contain a
complete set of representatives for ${\rm Log}(E_\fm)$ modulo
$p{\rm Log}(E_\fm)={\rm Log}(E_\fm^p)$.
Therefore, for any $M>0$ we may choose $\varepsilon_t\in E_\fm,\ t=1,\ldots,d-1$
such that ${\rm Log}(\varepsilon_t)\in B(l_t(M),r'(E_\fm))\ \forall\, t$
\emph{and} such that the matrix $A$ in $M_{d-1}(\bbZ)$ representing the $\varepsilon_t$'s
in the basis of the $\eta_t$'s has any required reduction $\bar{A}$
in $M_{d-1}(\bbZ/p\bbZ)$.
In particular we may insist that $\det(\bar{A})\neq 0$ so that
$p\ndiv |\det(A)|=|E_\fm\,:\,\langle\varepsilon_1,\ldots,\varepsilon_{d-1}\rangle|$.
By this means, we produce units $\varepsilon_1,\ldots,\varepsilon_{d-1}$
such that Condition\refeq{3X} holds with
$E:=\langle\varepsilon_1,\ldots,\varepsilon_{d-1}\rangle$. Moreover, if we substitute
$r'(E_\fm)$ for $r(E_\fm)$ throughout the proof of Lemme~2.1,
the reader may easily check that
its validity (with a new, corresponding choice of $M$) is unaffected.
Consequently the units $\varepsilon_1,\ldots,\varepsilon_{d-1}$ produced will
\emph{also} satisfy  Colmez' hypothesis $(H)$.
Now, for each $i\in\{1,\ldots,d\}$ and permutation $\sigma\in S_{d-1}$ Colmez defines
$f_{i,\sigma}:=\prod_{1\leq t<i}\varepsilon_{\sigma(t)}\in E$
(with $f_{1,\sigma}:=1\ \forall\,\sigma$). For any non-empty
$I\in\cP(\{1,\ldots,d\})$
(the power set) and $\sigma\in S_{d-1}$ we shall write $\bv_{(\sigma, I)}$
for the $|I|$-tuple consisting of the $f_{i,\sigma}$ ordered by increasing $i\in I$.
Then Lemme~2.2 of~\cite{Colm} says that a cone decomposition
$\tilde{\cC}$ for $E$ may be obtained by setting
$\tilde{\cC}:=\{\bv_{(\sigma, I)}\,:\,(\sigma, I)\in\cT\}$ where
$\cT$ is any set of representatives of equivalence classes for
$S_{d-1}\times \cP(\{1,\ldots,d\})$ modulo
a certain natural equivalence relation.
Finally, since $\ff'\neq \cO$, the character $\rho$ is
non-trivial, so there exists $b\in J$ with $\rho(b)\neq 1$.
Adding to $b$ a large positive element of $\ff' J\cap \bbZ$ if necessary,
we may assume that it also lies in $k^\times_\uinf$. Now
take $\cC$ to be the `multiplicative translate of $\tilde{\cC}$ by $b$'. In other words,
$\cC:=\{b\bv_{(\sigma, I)}\,:\,(\sigma, I)\in\cT\}$
where $b\bv_{(\sigma, I)}$ is the $|I|$-tuple consisting of the
$bf_{i,\sigma}$ ordered by increasing $i\in I$.
It is easy to see that $\cC$ is a cone decomposition
for $E$ (since $\tilde{\cC}$ is). Moreover, since $f_{i,\sigma}$ lies in
$E_\fm$, we have
$bf_{i,\sigma}\in J$ and
$bf_{i,\sigma}-b\in\ff J\subset\ff'J$ for each $i$ and $\sigma$.
Hence $\rho(bf_{i,\sigma})=\rho(b)\neq 1\ \forall\,i,\sigma$
which implies Condition\refeq{3Y} and thus completes the proof
of the Lemma.\ePf
\rem\label{rem:twobits}\

(i)~The above proof clarifies the purpose of
Conditions\refeq{3B} and\refeq{3C}: they allow
us to take $\ff'$ and hence $\rho$ to be non-trivial so that
the useful Condition\refeq{3Y} can be satisfied.

(ii)~In practice, the construction of Lemma~\ref{lemma:3L} may yield a subgroup $E$
whose index in $E_\fm$ is far larger than the minimum required to satisfy
Conditions\refeq{3X} and\refeq{3Y} for some $\cC$. Indeed,
in the interesting case $d=2$ ($k$ real quadratic) we can always take
$E$ to be $E_\fm$ itself (so, explicitly, $\cC=\{(b), (b,b\varepsilon)\}$ where
$E_\fm=\langle\varepsilon\rangle$ and $\rho(b)\neq 1$).
\vertsp\\
\noindent \textsc{Proof of Prop.~\ref{prop:3B} }\
Lemma~\ref{lemma:3L} gives a cone decomposition $\cC$ for $E\subset E_\fm$
such that
$
p\ndiv[E_\fm:E]
$
and \Hrhov\ holds for the given character $\rho$ and all $v\in\cC$.
We set
\[
F_\um(\uX):=\frac{1}{|E_\fm:E|}\sum_{\bv\in \cC}F_{\um,\bv}(\uX)\ \ \ \mbox{for each $\um\in\bbZ^d$.}
\]
Then
Condition~\ref{A'} of the Proposition follows from the first statement
of Lemma~\ref{lemma:3I.5} while
Condition~\ref{C'} follows from the second statement and the linearity of $V_{\ur}$.
Finally, Lemmas~\ref{lemma:3K} and~\ref{lemma:3J} imply that
$j(Z(0;\fw_{\rho,\us,\um}))$ lies in $\bbQ(\mu_{f'p^{n+1}})$ and that,
for any $\um\in\bbZ^d$ and $\us\in\cS(\un)$,
\begin{eqnarray*}
j(Z(0;\fw_{\rho,\us,\um}))&=&
\frac{1}{|E_\fm:E|}\sum_{\bv\in\cC}
j(Z_{\bv}(0,\xi_{\rho,\us,\um},\fp^\um J))\\
&=&\frac{1}{|E_\fm:E|}\sum_{\bv\in\cC}
F_{\um,\bv}(\zeta_n^{s_1}-1,\ldots,\zeta_n^{s_d}-1)
=F_{\um}(\zeta_n^{s_1}-1,\ldots,\zeta_n^{s_d}-1)
\end{eqnarray*}
which establishes
Condition~\ref{B'}. This completes the proof of Proposition~\ref{prop:3B}
hence also of Theorem~\ref{thm:3A}.\ePf
\section{Remarks and Conjectures}
\subsection{The $p$-Integrality of $\Phi_\fm(0)$}
Suppose that $\fm=\ff\uinf$ but $p$ does not necessarily split completely in $k$.
If $p\ndiv w_{k(\fm)}$ then \refeq{introA} implies that
$\Theta_{\fm(k(\fm))}(0)=\Theta_{k(\fm)/k}(0)$ lies in $\bbZ_p G$ and so
therefore must $\Theta_\fm(0)$ (which
is a $\bbZ G$-multiple of $\Theta_{\fm(k(\fm))}(0)$ by\refeq{2A}).
It follows from\refeq{2F.5} that
\beql{eq:4A}
j(\Phi_\fm(0))\in \bar{\bbZ}_p G
\eeq
whenever $p\ndiv w_{k(\fm)}$.
Suppose on the other hand that $p\neq 2$ splits completely in $k$ and that $\fm$ is as in\refeq{3C.5}
(in particular, the prime-to-$p$ part $\ff'$ of $\ff$ is non-trivial). Then
we may deduce from Lemma~\ref{lemma:3F}
and Propostion~\ref{prop:3B}
parts~\ref{A'} and~\ref{B'} that
$j(Z(0;\fw))$ lies in $\bar{\bbZ}_p$ for all $\fw\in\fW_\fm$ so\refeq{4A}
still holds even though we are assuming $p|\fm$, which is equivalent to $p|w_{k(\fm)}$
in this case. In fact, equation\refeq{2B} and Lemme~IV.1.1
of~\cite{Tate} can be used as in the proof of Theorem~\ref{thm:3A} to
reduce the condition on $\fm$ in this case to:\ \emph{either
$\ff'\neq \cO$ or $p\ndiv w_{k(\fm)}$}. The above considerations suggest the following
synthesis (with no splitting hypothesis on $p$).
\begin{conj} Suppose $k$ is totally real and $p\neq 2$. Then\refeq{4A} holds if
either $\ff'\neq \cO$ or $p\ndiv w_{k(\fm)}$.
\end{conj}
This conjecture might be established by adapting the techniques used
in the proof of~Theorem~\ref{thm:3A} to the case where $p$ does not split
completely in $k$. On the other hand, even when it does, \refeq{4A} will not
hold without \emph{some} condition on $\ff'$ and/or $w_{k(\fm)}$.
(Consider the case $k=\bbQ$, $\ff=\pnpo\bbZ$, \cf\ Example~\ref{ex:2A}).
The condition $p\neq 2$ may however be unnecessary.
\subsection{Theorem~\ref{thm:3A} Revisited}\label{subsec:thm3arev}
The question arises as to whether
the $p$-integrality property stated in Theorem~\ref{thm:3A} really is a new
phenomenon or whether it too might follow from a combination of\refeq{introA}
and Theorem~\ref{thm:2B} (presumably allowing us to dispose of the
splitting hypothesis). In support of the first alternative we now suppose that
$p$ splits in $k$ as $\fp_1\ldots\fp_d$ and examine
the case where $K=k(\fm)$ and $\fm=\ff\uinf$
is a conductor (for simplicity). In this case,
$G=G_\fm$, so Definition~\ref{def:sKkSKk},\refeq{2C} and
\refeq{2F.75} give (suppressing $j$ from the notation)
\[
\fs_{K/k}(\theta)=
\sum_{\fg|\ff}\frac{1}{\sqrt{d_k}N\fg}
  \left(\prod_{\fp|\ff,\ \fp\ndiv\fg}(1-N\fp\inv)\right)
  \nu_{\fn,\fm}(A_{\fn}\Theta_{\fn}(0))R_{K/k,p}(\theta)
\]
for any $\theta\in \bigwedge_{\bbZ G}^d U_p(K)$. (Recall that $\fn$ generically
denotes the cycle $\fg\uinf$).
Let us write $\Upsilon(\fm,\fg,\theta)$ for the term corresponding to $\fg$ in this
sum. It lies \emph{a priori} in $\barbbQ_p G$.
The most obvious hope is to use
\refeq{introA} to show that $\Upsilon(\fm,\fg,\theta)$ actually lies in
$\bar{\bbZ}_p G$ for all $\fg$ and thus deduce Theorem~\ref{thm:3A}.
We shall now show that in fact the
coefficients of $\Upsilon(\fm,\fg,\theta)$ may be arbitrarily large in $p$-adic
absolute value (and this for fixed $k$ and $\fg=\cO$). To do so, we fix also
an ideal $\fa$ of $\cO$ prime to $p$ such that $\fa\uinf$ is a conductor. Then $\fm_n:=p^{n+1}\fa\uinf$
is also a conductor for all $n\geq 0$ (that of $k(\fa\uinf)(\mu_{p^{n+1}})$). Write also
$K_n$ for $k(\fm_n)$ and $G_n$ for $\Gal(K_n/k)$. Then, taking $\fm=\fm_n$ and
$\theta_n\in\bigwedge_{\bbZ G_n}^d U^1_p(K_n)$, we have $A_\uinf=1$ and so
\begin{eqnarray*}
\Upsilon(\fm_n,\cO,\theta_n)&=&c\nu_{\uinf,\fm_n}(\Theta_{\uinf}(0))R_{K_n/k,p}(\theta_n)\\
                          &=&c\nu_{\uinf,\fm_n}\left(\Theta_{\uinf}(0)\pi_{\fm_n,\uinf}
                               (R_{K_n/k,p}(\theta_n))\right)\\
                          &=&c\nu_{\uinf,\fm_n}\left(\Theta_{\uinf}(0)
                              R_{k(\uinf)/k,p}(N_{k(\fm_n)/k(\uinf)}\theta_n)\right)\\
                            &=&c[k(\fm_n):k(\uinf)]\inv
                            \tilde{\nu}_{\uinf,\fm_n}\left(\Theta_{\uinf}(0)
                              R_{k(\uinf)/k,p}(N_{k(\fm_n)/k(\uinf)}\theta_n)\right)
\end{eqnarray*}
where $N_{k(\fm_n)/k(\uinf)}:\bigwedge_{\bbZ G_n}^d U^1_p(K_n)\rightarrow
\bigwedge_{\bbZ G_{\uinf}}^d U^1_p(k(\uinf))$ is the norm map
previously described and $c\in \bbQ_p^\times$ is independent of $n\geq 0$. The idea is now to
choose a sequence $\{\theta_n\}_{n\geq 0}$ such that
$N_{k(\fm_n)/k(\uinf)}\theta_n$ is a fixed element $\theta\in
\bigwedge_{\bbZ G_{\uinf}}^d U_p(k(\uinf))$ independent
of $n$ and such that $\Theta_{\uinf}(0)R_{k(\uinf)/k,p}(\theta)\neq 0$.
The desired result will then follow,
since $k(\fm_n)$ contains $k(\mu_{p^{n+1}})$ which implies that $p^n$ divides
$[k(\fm_n):k(\uinf)]$. The existence of such sequences is particularly
easy to establish under the following additional assumptions
\begin{enumerate}
\item\label{W} Leopoldt's Conjecture holds for $k$ at $p$,
\item\label{Y} $k(\uinf)$ is a totally complex quadratic extension of $k$ and
\item\label{Z} for each $i=1,\ldots,d$ there is a unique prime $\fP_i$ of $k(\uinf)$ lying
above $\fp_i$.
\end{enumerate}
These assumptions -- as well as the splitting hypothesis -- are valid
when $k=\bbQ(\sqrt{7})$ and $p=3$, for example. Conditions~\ref{Y} and~\ref{Z}
imply that $k(\uinf)_{\fP_i}/\bbQ_p$ is unramified of degree $2$ for each $i$
so that $U^1(k(\uinf)_{\fP_i})\cong \bbZ_p^2$ while Condition~\ref{W} implies
that the quotient of ${\displaystyle\lim_\leftarrow}\,G_n$ by a finite group
is isomorphic to $\bbZ_p$.
It follows from local class field theory that for each $i$ there exists
$u_i\in U^1(k(\uinf)_{\fP_i})\setminus\{ 1\}$ which is the
(local) norm of a principal local unit of $k(\fm_n)$ (completed at any prime above $\fP_i$)
for every $n\geq 0$.
Setting $\theta:=u_1\wedge\ldots\wedge u_d\in\bigwedge_{\bbZ G_{\uinf}}^d U^1_p(k(\uinf))$
we see that $\theta$ equals $N_{k(\fm_n)/k(\uinf)}\theta_n$
for some $\theta_n\in\bigwedge_{\bbZ G_n}^d U^1_p(k(\fm_n))$ for every $n\geq 0$.
It only remains to check that
$\Theta_{\uinf}(0)R_{k(\uinf)/k,p}(\theta)=\Theta_{\uinf}(0)
\lambda_{1,p}(u_1)\ldots\lambda_{d,p}(u_d)$ is non-zero.
Suppose this were not the case.
Condition~\ref{Y} implies that the unique  non-trivial character $\chi$ of $G_\uinf$
is totally odd so
$\chi(\Theta_{\uinf}(0))=L(0,\chi)\neq 0$ and we would have $\chi(\lambda_{i,p}(u_i))=0$
for some $i$.
On the other hand  $N_{k(\uinf)_{\fP_i}/k_{\fp_i}}u_i$ is, by construction,
a norm from $\bbQ(\mu_{p^{n+1}})$ (the completion of $k(\mu_{p^{n+1}})\subset k(\fm_n)$)
for every $n\geq 0$. Since $k_{\fp_i}\cong\bbQ_p$ it follows from local class field theory
that $N_{k(\uinf)_{\fP_i}/k_{\fp_i}}u_i=1$. Therefore the trivial character also
vanishes on $\lambda_{i,p}(u_i))$, so $\lambda_{i,p}(u_i)=0$
and in particular $\log_p(j\tau_i(u_i))=0$.
Since $\fP_i/p$ is unramified, this contradicts the fact
that $u_i\in U^1(k(\uinf)_{\fP_i})\setminus\{ 1\}$ and the desired result follows.

We note also a particular feature of our proof of Theorem~\ref{thm:3A}: By means
equations\refeq{3S.2} and\refeq{3S.3} it gives a neat and more-or-less
explicit expression for the coefficients
of $\fs_{K/k}(u_1\wedge\ldots\wedge u_d)$ where $u_i\in U^1(K_{\fP_i})$ for
$i=1,\ldots,d$. This assumes of course that the $g_i$'s and the $F_\ul$'s are known explicitly.
In practice, one might first choose a power series `$\tilde{h}$' as in
Lemma~\ref{lemma:3E} to determine each $u_i$ as well as the corresponding series
$g_i$. Furthermore, the proof of Proposition~\ref{prop:3B} would give a quite explicit
choice for $F_\ul$ once a suitable subgroup $E$ and cone decomposition $\cC$ had been found
(\eg\ by Remark~\ref{rem:twobits} for $d=2$).

Finally, on the basis of Theorem~\ref{thm:3A} and Corollary~\ref{cor:3A}
we hazard the
\begin{conj}\label{conj:4C} Suppose that $k$ is totally real and that $p$ is odd and splits
completely in $k$. Then
$\fS_{K/k}\subset p^{\delta(K/k)}\bbZ_p G$ for any $K$.
\end{conj}
Proposition~\ref{prop:unram} might be
taken as a (rather weak) hint that the splitting hypothesis can be dropped from
the above conjecture but
we prefer to keep it pending further evidence.
\subsection{Hilbert Symbols and Rubin-Stark Units}\label{subsec:rubinstark}
We assume from now on that the hypotheses of
Conjecture~\ref{conj:4C} hold and
that $K$ contains $\mu_{\pnpo}$ for some $n\geq 0$ so that
$K$ is totally complex and $\delta(K/k)=0$.
We shall also assume that the conclusions of
Conjecture~\ref{conj:4C} are valid (for instance, if\refeq{3B} holds, by Theorem~\ref{thm:3A}).
In this situation we shall set out some conjectural
congruences for $\fs_{K/k}$ modulo $\pnpo$ which
were foreshadowed in Example~\ref{ex:2A}
as generalisations of those coming from\refeq{secondfact} in the case $k=\bbQ$,
$K=\bbQ(\mu_\pnpo)$.

First, we need to define some pairings coming from local Hilbert symbols.
Let $K_p$ denote $\bbQ_p\otimes_\bbQ K=\prod_{\fP|p} K_\fP$ so that,
as in Section~\ref{sec:behav}, the composite
$j\tau_i$ extends to a map $K_p\rightarrow\barbbQ_p$ for
$i=1,\ldots,d$.
The field $\overline{j(\tau_i(K_p))}$ being
a finite (abelian) extension of $\bbQ_p$ containing $\mu_\pnpo$, the
Hilbert symbol
$(\cdot,\cdot)_\pnpo$ gives a well-defined, $\mu_\pnpo$-valued,
bilinear pairing on its multiplicative group.
Thus for each $i$ we obtain a unique $\bbZ$-bilinear pairing
$[\cdot,\cdot]_{i,n}: K_p^\times \times K_p^\times\rightarrow \bbZ/\pnpo\bbZ$ satisfying
\[
\zeta_n^{[\alpha,\beta]_{i,n}}=(j\tau_i(\alpha),j\tau_i(\beta))_\pnpo
\ \ \ \mbox{for all $\alpha,\beta\in K_p^\times$}
\]
This exhibits the following equivariance property with respect to
the decomposition group $D_{\fp_i}=G(K/k)^{-1}_{\fp_i}\subset G$
\beql{eq:4B}
[\delta\alpha,\delta\beta]_{i,n}=\kappa_n(\delta)[\alpha,\beta]_{i,n}\ \ \ \mbox{for all
$\delta\in D_{\fp_i}$}
\eeq
where $\kappa_n:G\rightarrow(\bbZ/\pnpo\bbZ)^\times$ is the character defined by
$\sigma(\zeta)=\zeta_n^{\kappa_n(\sigma)}\ \forall\,\zeta\in\mu_\pnpo,\sigma\in G$.
Choose a set $R_i$ of coset representatives for $D_{\fp_i}$ in $G$ and
define a pairing
$[\cdot,\cdot]^G_{i,n}$ on $K_p^\times$
with values in the group-ring $(\bbZ/\pnpo\bbZ)G$ by setting
\[
[\alpha,\beta]^G_{i,n}:=
\sum_{\rho\in R_i}
\kappa_n(\rho)\inv\sum_{\sigma\in G}[\rho\alpha,\sigma\beta]_{i,n}\sigma\inv\rho
\ \ \ \mbox{for all $\alpha,\beta\in K_p^\times$}
\]
It follows from property\refeq{4B} that $[\cdot,\cdot]^G_{i,n}$
is independent of the choice of $R_i$. Moreover it is
$\bbZ G$-(semi)linear in each variable as follows.
For any $\alpha,\beta\in K_p^\times$ we have
\beql{eq:4C}
[x\alpha,\beta]^G_{i,n}=
\kappa_n^\ast(x)[\alpha,\beta]^G_{i,n}
\ \ \mbox{and}\ \
[\alpha,x\beta]^G_{i,n}=
x[\alpha,\beta]^G_{i,n}\ \ \ \mbox{for all $x\in \bbZ G$}
\eeq
where $\kappa_n^\ast:\bbZ G\rightarrow (\bbZ/\pnpo \bbZ)G$ is the ring homomorphism
sending $\sum a_g g$ to $\sum \overline{a_g}\kappa_n(g)g\inv$.
(We shall summarise property\refeq{4C} by saying that the pairing
$[\cdot,\cdot]^G_{i,n}$ is \emph{$(\kappa_n^\ast,1)$-bilinear}.)
Summing over $i$ we get a
pairing $[\cdot,\cdot]^G_n:K_p^\times\times K_p^\times\rightarrow(\bbZ/\pnpo\bbZ)G$,
namely $[\alpha,\beta]^G_n:=\sum_{i=1}^d[\alpha,\beta]^G_{i,n}$, which
is also $(\kappa_n^\ast,1)$-bilinear.
Finally, by regarding both $K^\times$ and $U^1_p(K)$
a submodules of $K_p^\times$, it follows that there is a
unique pairing
$\cH_n$ from
$\bigwedge^d_{\bbZ G}K^\times\times\bigwedge^d_{\bbZ_p G} U^1_p(K)$
to $(\bbZ/\pnpo\bbZ)G$ which
is $(\kappa_n^\ast,1)$-bilinear
and satisfies
\[
\cH_n(x_1\wedge\ldots\wedge x_d,u_1\wedge\ldots\wedge u_d)=\det([x_i,u_t]^G_n)_{i,t}
\]
We would now like to construct an element $\eta\in\bigwedge^d_{\bbZ G}K^\times$
such that $\cH_n(\eta,\cdot)$ plays essentially the same r\^ole as
$((1-\zeta_n)\inv,\cdot)_\pnpo$ plaed in equation\refeq{secondfact} of
Example~\ref{ex:2A} (where $d=1$).
The best tool we currently have for this construction is
Rubin's restatement of the Stark Conjecture for $K^+/k$
where $K^+$ is the maximal real subfield of $K$.
Not only is this conjecture itself unproven in almost all the relevant cases
but, as we shall see, its
very nature makes our conjectural congruences for $\sKk$ vaguer
and more awkward than Example~\ref{ex:2A} might
reasonably lead one to expect.
With these \emph{caveats}, let us recall Rubin's formulation using
the higher derivatives of the function $\Theta_{K^+/k}(s)$ at $s=0$,
referring to~\cite{Ru} and~\cite[\S4]{twizas} for more details and relevant remarks.
(See also~\cite[\S5]{twizas} and~\cite{sconzp} for a reformulation using
$\Phi_{K^+/k}(s)$ at $s=1$, essentially by Theorem~\ref{thm:2A}).
Let us write $S_0$ for $S_{\rm ram}(K/k)$ and
$U_{S_0}(K^+)$ for the group of (global) $S_0$-units of $K^+$.
Our assumptions force $\fp_i\in S_0\ \forall\,i$. This means that
$|S_0|\geq 2d$ (and also that $U_{S_0}(K^+)$
is \emph{not} contained in the product of local units, $U_p(K^+)$, despite the notation).
Set $G^+:=\Gal(K^+/k)$ and define the archimedean analogue of $\lambda_{K^+/k,i,p}$ by
$\lambda_{K^+/k,i}(a)=\sum_{g\in G^+}\log|\tau_i(ga)|g\inv\in \bbR G^+$
for any $a\in K^{+,\times}$. We obtain a
unique, archimedean regulator
$R_{K^+/k}:\bbQ\otimes\bigwedge^d_{\bbZ G^+}U_{S_0}(K^+)\rightarrow \bbR G^+$
taking $u_1\wedge\ldots\wedge u_d$ to $\det(\lambda_{K^+/k,i}(u_t))_{i,t=1}^d$.
On the other hand, we may define
\[
\Theta_{K^+/k, S_0}(s):=
\prod_{\fp\in S_0\atop \fp\nin S_{\rm ram}(K^+/k)}
(1-N\fp^{-s}\sigma_{\fp,K^+}\inv)\Theta_{K^+/k}(s)=
\pi_{K,K^+}(\Theta_{K/k}(s))
\]
(This is the function $\Theta_{K^+,S_0,\emptyset}(s)$ of~\cite{Ru}).
Since $S_0$ contains at least $d+1$ places of $k$,
and the the $d$ real ones split completely in $K^+$,
it follows that $\Theta_{K^+/k, S_0}(s)$ has at
least a $d$-fold zero at $s=0$. We set
$\Theta_{K^+/k, S_0}^{(d)}(0):=\lim_{s\rightarrow 0}s^{-d}\Theta_{K^+/k, S_0}(s)\in \bbC G^+$
and write $X(S_0,d,G^+)$ for the set of complex
irreducible characters $\chi$ of $G^+$ for which
$\chi(\Theta_{K^+/k, S_0}^{(d)}(0))\neq 0$. Let us take the extension
`$K/k$' of~\cite{Ru} to be $K^+/k$,`$S$' to be $S_0$,`$T$' to be $\emptyset$, `$r$' to be $d$
and the chosen places
`$w_1,\ldots,w_d$' of $K^+$ to be the real ones defined by $\tau_1,\ldots,\tau_d$.
The torsion subgroup of $U_{S_0}$ is $\{\pm 1\}$ so the conditions
of Rubin's Conjectures B and B$'$ `over $\bbZ$' (see~\cite[Hyp.~2.1]{Ru})
are actually not quite met with $T=\emptyset$.
However, we shall only need his Conjecture~A$'$ `over $\bbQ$',
for which the choice of $T$ is irrelevant. The latter is in fact equivalent to
Stark's conjecture for each
$\chi\in X(S_0,d,G^+)$ (see~\cite[Conjecture I.5.1]{Tate}) and states
\begin{conj}[Conjecture A$'$ of~\cite{Ru}]\label{conj:Rubin}\ \\
There exists an element $\eta_{K^+/k,S_0}$
of $\bbQ\otimes\bigwedge^d_{\bbZ G^+}U_{S_0}(K^+)$ such that
\begin{enumerate}
\item
$\Theta_{K^+/k, S_0}^{(d)}(0)=R_{K^+/k}(\eta_{K^+/k,S_0})$
and
\item
$e_\chi\eta_{K^+/k,S_0}=0$ in $\bbC\otimes\bigwedge^d_{\bbZ G^+}U_{S_0}(K^+)$
for every character $\chi$ of $G^+$ not in $X(S_0,d,G^+)$.
\end{enumerate}
\end{conj}
We assume henceforth that $\eta_{K^+/k,S_0}$ satisfies the two conditions of this
conjecture. This actually makes $\eta_{K^+/k,S_0}$ \emph{unique}
and we define $\eta^+_{K/k}$ to be its image under the natural map
$\bbQ\otimes\bigwedge^d_{\bbZ G^+}U_{S_0}(K^+)\rightarrow
\bbQ\otimes\bigwedge^d_{\bbZ G}K^\times$.
Consider the case $K/k=\bbQ(\zeta_\pnpo)/\bbQ$.
We have $d=1$, $S_0=\{\infty,p\}$ and
one can check that Conjecture~\ref{conj:Rubin} always holds with
$\eta_{K^+/\bbQ,S_0}=-\half\otimes(1-\zeta_n)(1-\zeta_n\inv)$ (\cf\ \cite[\S III.5]{Tate} and
\cite[\S 3.5]{zetap1}.
In fact, $\Theta_{K^+/\bbQ,S_0}(s)=\Theta_{K^+/\bbQ}(s)$ except in the
trivial case $\pnpo=3$). It follows that $\eta^+_{K/\bbQ}=1\otimes(1-\zeta_n)\inv
\in\bbQ\otimes\bigwedge^1_{\bbZ G} K^\times=\bbQ\otimes K^\times$.
Moreover the reader may check
that if we take $\tau_1$ to be $1\in\Gal(\barbbQ/\bbQ)$ then\refeq{secondfact}
is precisely the statement that
$\overline{\fs_{K/\bbQ}}=\cH_n((1-\zeta_n)\inv,u)$
in $(\bbZ/\pnpo\bbZ)G$ for all $u\in U^1_p(K)=\bigwedge^1_{\bbZ_p G} U^1_p(K)$.
Unfortunately, for general $K^+/k$, Rubin's conjecture
does not require that $\eta^+_{K/k}$ should be
of the form $1\otimes\tilde{\eta}$ for some $\tilde{\eta}\in\bigwedge^d_{\bbZ G} K^\times$.
The author knows of no case satisfying our hypotheses where this condition
actually fails. On the other hand, \cite[Prop.~4.4]{Ru} shows that in a very different
case, the counterpart of the element $\eta_{K^+/k,S_0}$
satisfying Rubin's general conjecture need not even lie in the counterpart of
$\bbZ[\half]\otimes\bigwedge^d_{\bbZ G^+}U_{S_0}(K^+)$.  It therefore seems safer to
proceed as follows.
We define an ideal $\cI(\eta^+_{K/k})$ of finite index in $\bbZ G$ by
\[
\cI(\eta^+_{K/k}):=\{x\in \bbZ G :
\mbox{
$x\eta^+_{K/k}=
1\otimes\tilde{\eta}$
for some
$\tilde{\eta}\in\bigwedge^d_{\bbZ G} K^\times$
\}
}
\]
and formulate the
\begin{conj}\label{conj:4D} Suppose that $p$ is odd and splits completely in $k$
and that other hypotheses and notations are as above. Then for
every $x\in\cI(\eta^+_{K/k})$ there exists
$\tilde{\eta}_x\in\bigwedge^d_{\bbZ G}K^\times$ such that
\begin{enumerate}
\item\label{part:conj4D1}
$x\eta^+_{K/k}=1\otimes\tilde{\eta}_x$ and
\item\label{part:conj4D2}
$\kappa_n^\ast(x)\overline{\sKk(\theta)}=\cH_n(\tilde{\eta}_x,\theta)$
in $(\bbZ/\pnpo\bbZ)G$ for all
$\theta\in\bigwedge^d_{\bbZ_p G} U^1_p(K)$.
\end{enumerate}
\end{conj}
It clearly suffices to check this conjecture for a set of elements $x$ generating
$\cI(\eta^+_{K/k})$ over $\bbZ G$.
Note that condition~\ref{part:conj4D1} only
determines $\tilde{\eta}_x$ up to
$\bbZ$-torsion in $\bigwedge^d_{\bbZ G}K^\times$, which
does not lie in the kernel of $\cH_n(\cdot,\theta)$ for all $\theta$, even in
the case $K/k=\bbQ(\zeta_\pnpo)/\bbQ$. This means
firstly that condition~\ref{part:conj4D2} cannot be expected to
hold \emph{for all lifts $\tilde{\eta}_x$} satisfying
condition~\ref{part:conj4D1} and secondly that one could
weaken the conjecture by allowing $\tilde{\eta}_x$ satisfying the equality
in~\ref{part:conj4D2} to depend on $\theta$ as well as $x$
(subject to~\ref{part:conj4D1}).

We briefly explain the two pieces of evidence that
motivate Conjecture~\ref{conj:4D}, postponing a
more detailed account for a future paper. Firstly, if $k=\bbQ$
then the conjecture holds with $\cI(\eta^+_{K/k})=\bbZ G$ and $\tilde{\eta}_1=
N_{\bbQ(\xi_{f(K)})/K}(1-\xi_{f(K)})\inv$ where $\xi_{f(K)}=\be(1/f(K))$.
This may be deduced from equation\refeq{2K} and a generalisation of Artin-Hasse's
explicit reciprocity law to $\bbQ_p(\mu_{f(K)})$ which is proved by Coleman
in~\cite{Coleman2}.
Secondly, Conjecture~\ref{conj:4D} actually follows from
the same result of Coleman together with those of~\cite{sconzp}
\emph{provided} that it is weakened by using a slightly smaller ideal
than $\cI(\eta^+_{K/k})$ \emph{and}
that certain extra hypotheses are satisfied. The most important of these
are: (a)~that $K=L(\mu_\pnpo)$ for some totally real or CM extension
$L$ of $k$ which is unramified above $p$ and (b)~that a certain
strongly norm-coherent form of Conjecture~\ref{conj:Rubin} over
$k$, together with its $p$-adic analogue, holds in the cyclotomic $\bbZ_p$-tower over $K^+$.
For a more precise form of hypothesis~(b) we refer to properties
$P1(K^+/k,p)$, and \eg\ $P2(K^+/k,p)$ of~\cite[\S 3]{sconzp}. By using Example~3.2
of~\emph{ibid.}, one can construct a large infinite class of
extensions $K/k$ with $K$ absolutely abelian and $k\neq \bbQ$ which satisfy all the above
hypotheses and hence also the weakened form of Conjecture~\ref{conj:4D}.

Finally, we mention that methods from \eg~\cite{zetap2} could be used to
test Conjecture~\ref{conj:4D} numerically and also
the necessity of the condition that $p$ split completely in $k$.

\end{document}